\documentclass[9pt]{amsart}

\usepackage{amsmath, amssymb}
\usepackage{amssymb}
\usepackage{amscd}
\usepackage{amsthm}
\usepackage{eucal}
\usepackage{graphicx}
\usepackage{amsfonts}
\usepackage{latexsym}
\numberwithin{equation}{section}

\def\s{\sigma}
\def\a{\alpha}
\def\b{\beta}
\newtheorem{thm}{Theorem}

\newtheorem{prop}{Proposition}[section]

\newtheorem{lem}{Lemma}[section]
\newtheorem{defi}{Definition}[section]

\newtheorem{cor}{Corollary}
\newtheorem{example}{Example}

\newtheorem*{remark}{Remark}
\begin{document}
\title{Partially Ordering Unknotting Operations}
\author{Maki Nagura}
\address{Yokohama National University, Tokiwadai, Hodogaya-ku, Yokohama 240-8501, Japan}
\email{maki{@}ynu.ac.jp}
%\subjclass[2010]{55R40, 57R67}
\maketitle

\begin{abstract}
In this paper, we introduce an equivalence relation on the set of local moves and classify local moves, called the extended $ST$-moves, up to the equivalence. Moreover, by inducing a binary relation on the set of equivalence classes of local moves, we show that an extended $ST$-move realizes the crossing change or the $SH(2)$-move. In addition, for any oriented knot and two extended $ST$-moves, we disscus the magnitude relation between the unknotting numbers of the knot via the moves, and show that there is an extended $ST$-move except $SH$-moves so that the knot can be transformed into the trivial knot by the single extended $ST$-move. Finally, we provide some examples of $ST$-moves with the binary relation.

\end{abstract}

\section{introduction}

An operation that replaces a tangle diagram on a knot or link diagram with another tangle diagram is referred to as ``a local move on a knot or link diagram." For instance, the Reidemeister moves \cite{re} are local moves on a knot or link diagram. The crossing change, called the $X$-move, the $\Delta$-move \cite{MN}, the $\Delta_{ij}$-move \cite{N}, the $\sharp$-move \cite{mura} and the $n$-gon move \cite{miya} are also local moves on a knot or link diagram. We will define local moves as pairs of two tangle diagrams (see Definition 2.5).

\begin{figure}[h]
\begin{center}
\includegraphics[scale=0.7]{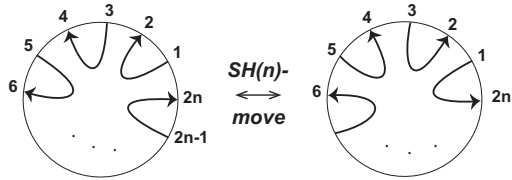}
\end{center}
\caption{$SH(n)$-move.}
\label{shn}
\end{figure}

\begin{figure}[h]
\begin{center}
\includegraphics[scale=0.7]{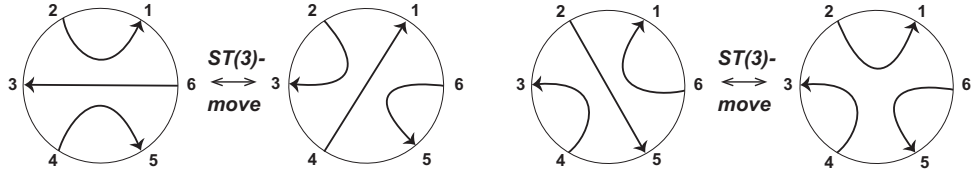}
\end{center}
\caption{Examples of $ST(3)$-moves.}
\label{st3new}
\end{figure}

In \cite{HNT}, J. Hoste, Y. Nakanishi and K. Taniyama defined an $SH(n)$-move (see FIG. \ref{shn}), and verified that an $SH(2n-1)$-move is an unknotting operation on an oriented knot or link diagram for $2\leq n\in\mathbb{N}$. An unknotting operation (see, e.g., \cite{kawauchi}) is a local move on a knot or link diagram such that any knot or link diagram can be transformed into a trivial knot or link diagram by a finite sequence of the local move and Reidemeister moves. In \cite{na}, we defined an $ST(n)$-move, which is an extension of an $SH(n)$-move, and demonstrated that it realizes the crossing change or the $SH(2)$-move for $2\leq n\in\mathbb{N}$. A local move is called an $ST(n)$-move if the two oriented $n$-tangle diagrams are both trivial and not equal (see \cite{na} and FIG. \ref{st3new}).

In this paper, we introduce an equivalence relation on the set of local moves and classify local moves, called the extended $ST$-moves (see Definition \ref{def2-8}), up to the equivalence. In Theorem 1, we prove that there is a one-to-one correspondence between the set of equivalence classes of extended $ST$-moves and the set of standard $ST$-moves (see Definition \ref{def3-3}). Therefore, any standard $ST$-move can be choosed as a representative of an equivalence class of extended $ST$-moves. In Theorem \ref{thm2}, we show that an extended $ST$-move realizes the crossing change or the $SH(2)$-move. The former extended $ST$-move is called an $X$-type and the latter one is called an $O$-type. In Theorem 3, it is shown that a local move, which realizes an $X$-type, is an unknotting operation. Given any knot and two $X$-types, we obtain the magnitude relation between the unknotting numbers of the knot via the two moves (Theorem \ref{thm4}). Moreover, we show that for any oriented knot $K$, there is an extended $ST$-move except $SH$-moves so that any diagram of $K$ can be transformed into a trivial knot diagram by the single $ST$-move (Theorem \ref{thm5}).

Section 2 presents some definitions and a proposition necessary for proving Theorem 1. In section 3, we describe some lemmas and prove Theorem 1. In section 4, we introduce a binary relation $\preceq$ on the set of equivalence classes of local moves and demonstrate that the binary relation is a partial order on the set. In section 5, we state necessary and sufficient conditions (Lemmas \ref{lem5-2} and \ref{lem5-3}) for the partial order between the equivalence class of an $ST$-move and one of an $SH$-move to exist by using their representatives. Section 6 discusses the unknotting numbers of $ST$-moves. Finally, we provide some examples of $ST$-moves with the relation $\preceq$.

\section{Definitions}

Throughout this paper, we work in PL category. Tangles were introduced by J. Conway in \cite{co} in order to help in assembling a knot table and develop symbols of knot diagrams. Since then, tangles have been useful for studying knot theory, DNA topology, quantum topology, and applied fields, such as molecular biology. 

We shall begin with some definitions on tangles. First, we will define a tangle as follows.

\begin{defi}\label{def2-1}
Let $B=\{(x,y,z)\in\mathbb{R}^3 \,|\, x^2+y^2+z^2\leq 1\}$ be a unit 3-ball. Let $\Hat{T}=\cup_{i=1}^n t_i$ be a union of $n$ pairwise disjoint arcs $t_i$ embedded properly in $B$ and let $\partial \Hat{T}=\partial(\cup_{i=1}^n t_i)=\cup_{i=1}^n\partial t_i=\bigl\{(\cos{ j\over n}\pi , \, \sin {j\over n}\pi , \, 0)\in\mathbb{R}^3\, |\, j=1,2,\ldots , 2n\bigr\}$. Then $(B,\Hat{T})$ is called an $n$-tangle. An $n$-tangle $(B, \Hat{T})$ is called to be oriented if each arc $t_i$ is oriented, where $i=1,2,\ldots , n$.
\end{defi}

In this paper, we treat the following regular diagrams.

\begin{defi} \label{def2-2} Let $(B, \Hat{T})$ be an oriented $n$-tangle. Let $p$ be a projection of $B$ onto the unit disk $D=\{(x,y,0)\,|\, x^2+y^2\leq 1\}.$ Then $(D, T):=(p(B), p(\Hat{T}))$ is called a tangle diagram of $(B, \Hat{T})$, where $p(\Hat{T})$ is a regular diagram of $\Hat{T}$, i.e., $p(\Hat{T})$ is regular, and we draw one arc close to a double point (or crossing) so that it appears to have been cut to express that the arc passes under the other arc. Each point $\bigl(\cos{ j\over n}\pi , \, \sin { j\over n}\pi ,\, 0\bigr)$ is marked $j$ and called an e-point of $(B, \Hat{T})$ or $(D, T)$.\end{defi}

Henceforth, we let $p$ be a projection of $B$ onto the unit disk $D=\{(x,y,0)\,|\, x^2+y^2\leq 1\}$, and assume that tangle diagrams are oriented. In the next definition, we will describe the equality of two tangle diagrams. 

\begin{defi} \label{def2-3}
Let $(D_1, T_1)$ and $(D_2, T_2)$ be $n$-tangle diagrams and let $I(\partial T_1)=I(\partial T_2)$, where $I(\partial T_i)$ is the set of initial points of $\partial T_i$. If we can change $T_1$ into $T_2$ by performing a finite number of Reidemeister moves in $D_1$ and $D_2$, respectively, keeping the $2n$ marked e-points fixed, then the tangle diagrams $(D_1, T_1)$ and $(D_2, T_2)$ are also said to be equal and are denoted by $(D_1, T_1)=(D_2, T_2)$ or $T_1=T_2$. 
\end{defi}

\begin{defi} \label{def2-4}
A tangle diagram $(D,T)$ is trivial if we can change $(D, T)$ into a diagram with no crossings by performing a finite number of Reidemeister moves in $D$, keeping the marked e-points fixed.
\end{defi}

\begin{defi} \label{def2-5}
A local move is a pair of tangle diagrams, $(D_1, T_1)$ and $(D_2, T_2)$, with $\partial T_1=\partial T_2$ and $I(\partial T_1)=I(\partial T_2)$. It is denoted by $\mathcal{L}: (D_1, T_1)\leftrightarrow (D_2, T_2)$, $\mathcal{L}: T_1\leftrightarrow T_2$ or simply denoted by $\mathcal{L}$. The set of local moves is denoted by $\mathbb{L}$.
\end{defi}

Let $(D, T)$ be a tangle diagram and $f$ be a map from $D$ to $\widetilde{D}=\{(x,y,0)\in\mathbb{R}^3 \,|\, x^2+y^2\leq 1/2\}$ such that $f((x,y,0))={1\over 2}(x,y,0)$. Then we let $(\widetilde{D}, \widetilde{T}):=(f(D), f(T))$ and call $(\widetilde{D}, \widetilde{T})$ the shrunk diagram of $(D, T)$. 

In order to describe the equivalence of local moves, we shall define a specific operation, called a braiding operation, as follows.

\begin{figure}[h]
\begin{center}
\includegraphics[scale=0.5]{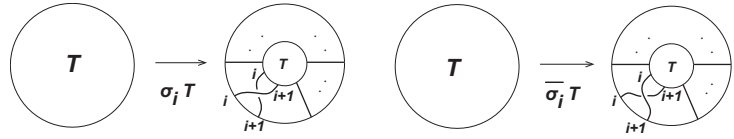}
\end{center}
\caption{Braiding operations $\s_i$ and $\overline{\s_i}.$}
\label{eq-local}
\end{figure}

\begin{defi} \label{def2-6}
Let $i\, (\leq 2n)$ be a positive integer and $A:=\{(x,y,0)\in\mathbb{R}^3 \,|\, 1/2\leq x^2+y^2\leq 1\}$. Let $\Hat{E}_i=\cup_{k=1}^{2n} {e}_k$ and $\Hat{E}_{\overline{i}}=\cup_{k=1}^{2n} \overline{{e}_k}$ be unions of $2n$ pairwise disjoint arcs $e_k$ and $\overline{{e}_k}$ embedded properly in $B\setminus Int(\widetilde{B})$, respectively, satisfying the following:

\noindent $(I)$ Case of $i=1,2,\ldots , 2n-1,$
\begin{itemize}
\item $e_j$ and $\overline{{e}_j}$ are both the lines connecting the points $1/2(\cos{ j\over n}\pi , \, \sin {j\over n}\pi , \, 0)$ and $(\cos{ j\over n}\pi , \, \sin { j\over n}\pi , \, 0)$ for $j\neq i,i+1, $

\item $e_i$ and $\overline{{e}_i}$ are arcs connecting the points $1/2(\cos{ i\over n}\pi , \, \sin { i\over n}\pi , \, 0)$ and $(\cos{i+1\over n}\pi ,$ $\sin {i+1\over n}\pi$, $\, 0)$ whose images $p(e_i)$ and $p(\overline{{e}_{i}})$ are on the annulus $A$, 

\item $e_{i+1}$ and $\,\overline{{e}_{i+1}}$ are arcs connecting the points $1/2(\cos{i+1\over n}\pi ,\,  \sin {i+1\over n}\pi , \, 0)$ and $($$\cos{i\over n}\pi ,$ $\, \sin {i\over n}\pi ,$ $\, 0)$ whose images $p(e_{i+1})$ and $p(\overline{{e}_{i+1}})$ are on the annulus $A$, 

\item the diagram $(A, E_i:=p(\Hat{E}_i))$ has only one crossing which is an overcrossing $($or an undercrossing$)$ on $p(e_i)$ $($or $p(e_{i+1}))$ and 
\item the diagram $(A, E_{\overline{i}}:=p(\Hat{E}_{\overline{i}}))$ has only one crossing which is an undercrossing $($or an overcrossing$)$ on $p(\overline{e_i})$ $($or $p(\overline{e_{i+1}}))$. 
\end{itemize}

\noindent $(II)$ Case of $i=2n,$
\begin{itemize}
\item $e_j$ and $\overline{e_j}$ are both the lines connecting the points $1/2(\cos{ j\over n}\pi , \, \sin {j\over n}\pi , \, 0)$ and $(\cos{ j\over n}\pi , \, \sin { j\over n}\pi , \, 0)$ for $j\neq 1,2n, $
\item $e_{2n}$ and $\overline{e_{2n}}$ are arcs connecting the points $1/2(1, \, 0 , \, 0)$ and $(\cos{1\over n}\pi , \, \sin {1\over n}\pi , \, 0)$ whose image $p(e_{2n})$ and $p(\overline{{e}_{2n}})$ are on $A$, 
\item $e_1$ and $\overline{e_1}$ are arcs connecting the points $1/2(\cos{1\over n}\pi ,\,  \sin {1\over n}\pi , \, 0)$ and $(1 , \, 0 , \, 0)$ whose images $p(e_{1})$ and $p(\overline{{e}_{1}})$ are on $A$,  
\item the diagram $(A, E_{2n}:=p(\Hat{E}_{2n}))$ has only one crossing which is the overcrossing $($or undercrossing$)$ on the arc $p(e_{2n})$ $($or $p(e_{1})$$)$ and
\item the diagram $(A, E_{\overline{2n}}:=p(\Hat{E}_{\overline{2n}}))$ has only one crossing which is the undercrossing $($or overcrossing$)$ on the arc $p(\overline{e_{2n}})$ $($or $p(\overline{e_{1}}))$. 
\end{itemize}
If $(D, T)$ is an $n$-tangle diagram and $(\widetilde{D}, \widetilde{T})$ is the shrunk diagram of $(D, T)$, then we can regard $(D_1, T_1):=(A\cup\widetilde{D},\, E_i\cup\widetilde{T})$ and $(D_2, T_2):=(A\cup\widetilde{D},\, E_{\overline{i}}\cup\widetilde{T})$ as $n$-tangle diagrams. The operation that transforms $(D, T)$ into $(D_1, T_1)$ is called the braiding operation $\s_i$ on $(D,T)$ and we write $T_1=\s_i\, T$. Similarly, the operation that transforms $(D, T)$ into $(D_2, T_2)$ is called the braiding operation $\overline{\s_i}$ on $(D,T)$ and we write $T_2=\overline{\s_i}\, T$ $($see FIG. \ref{eq-local}$)$. Here the orientation of $(D_j, T_j)$, $j=1,2$, is induced from the orientation of $(D, T)$.

\end{defi}

From the definition, we see that $\s_2(\s_1 T)=(\s_2\s_1)T$. Therefore, we denote $\s_2(\s_1 T)$ by $\s_2\s_1 T$.

\begin{defi}\label{def2-7}
Two local moves, $\mathcal{L}: T_1\leftrightarrow T_2$ and $\mathcal{L}':T_1'\leftrightarrow T_2'$, are equivalent, denoted by $\mathcal{L}\cong\mathcal{L}'$, if there exists a finite sequence of braiding operations $\s_{i_1}$, $\s_{i_2}$, $\cdots$, $\s_{i_m}$ such that $T_1'=\s_{i_m}\cdots \s_{i_2}\s_{i_1}T_1$ and $T_2'=\s_{i_m}\cdots \s_{i_2}\s_{i_1}T_2$. Then the operation that transforms $\mathcal{L}$ into $\mathcal{L}'$ is called a sequence of braiding operations $\s_{i_1}$, $\s_{i_2}$, $\cdots$, $\s_{i_m}$ on $\mathcal{L}$ and we say that $\mathcal{L}'$ can be obtained from $\mathcal{L}$ using a finite sequence of braiding operations $\s_{i_1}$, $\s_{i_2}$, $\cdots$, $\s_{i_m}$.
\end{defi}

This relation $\cong$ is clearly an equivalence relation, i.e., it satisfies the following properties $(i)-(iii)$: $(i)$ $\mathcal{L}\cong\mathcal{L}$, $(ii)$ $\mathcal{L}\cong\mathcal{L}'$ implies $\mathcal{L}'\cong\mathcal{L}$ and $(iii)$ $\mathcal{L}\cong\mathcal{L}'$ and $\mathcal{L}'\cong\mathcal{L}''$ imply $\mathcal{L}\cong\mathcal{L}''$. If $T_1'=\s_{i_m}\cdots \s_{i_2}\s_{i_1}T_1$ and $T_2'=\s_{i_m}\cdots \s_{i_2}\s_{i_1}T_2$, then $T_1=\overline{\s_{i_1}}\cdots \overline{\s_{i_{m-1}}}\, \overline{\s_{i_m}}\, T_1'$ and $T_2=\overline{\s_{i_1}}\cdots \overline{\s_{i_{m-1}}}\, \overline{\s_{i_m}}\, T_2'$. 

\begin{figure}[h]
\begin{center}
\includegraphics[scale=0.5]{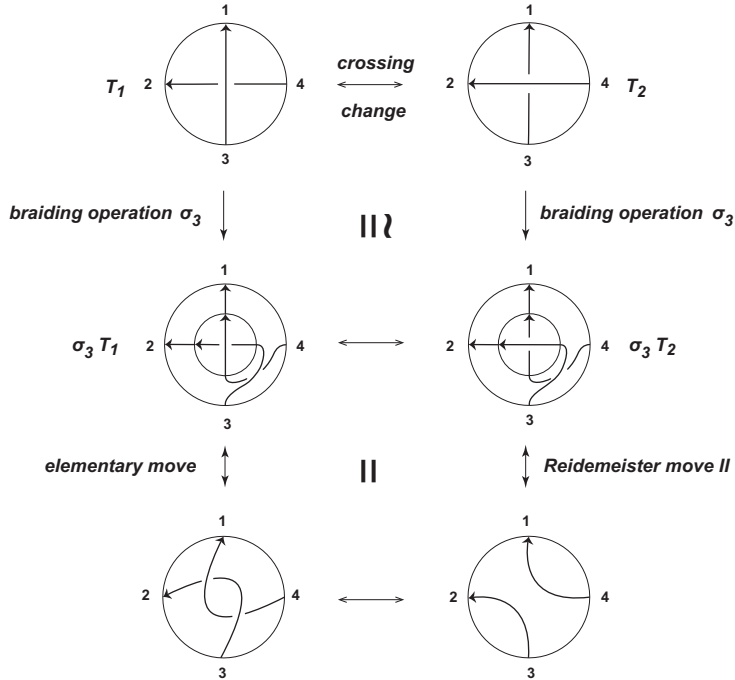}
\end{center}
\caption{Equivalent local moves.}
\label{crossch}
\end{figure}

 For example, the crossing change ($X$-move) is equivalent to the local move, as shown in the center-side or lower-side diagram of FIG. \ref{crossch}. Even if a rotation is performed on two tangle diagrams of a local move, the (old) local move and the new local move are equivalent. 

The next proposition demonstrates that a local move and the local move after a rotation are equivalent.

\begin{prop}\label{prop2-1}
Let $\mathcal{L}:T_1\leftrightarrow T_2$ be a local move and let $T'_i$ be the $n$-tangle diagram rotated $T_i$ ${m\over n}\pi$ radians around its origin $(m\in\mathbb{Z}$, $i=1, 2)$. Then $\mathcal{L}':T'_1\leftrightarrow T'_2$ is a local move such that $\mathcal{L}\cong\mathcal{L}'$. 
\end{prop}

\begin{proof} In order to prove the proposition, it is helpful to separate our proof into three cases: $m=0$, $m>0$ and $m<0$. 

Case $m=0$: Then $T'_i=T_i$. Thus, $\mathcal{L}=\mathcal{L}'$.

Case $m>0$: Let $T''_i=\s_1\s_2\cdots \s_{2n-1}\, T_i$. Then the tangle diagram $T''_i$ is one rotated $T_i$, ${1\over n}\pi$ radians around its origin. Thus, $T'_i=(\s_1\s_2\cdots \s_{2n-1})^m\, T_i$. Thus, we have $\mathcal{L}\cong\mathcal{L}'$.

Case $m<0$: Let $T''_i=\overline{\s_{2n-1}}\, \overline{\s_{2n-2}}\cdots \overline{\s_1}\, T_i$. Then the diagram $T''_i$ is one rotated $T_i$, $-{1\over n}\pi$ radians around its origin. Thus, $T'_i$$=(\overline{\s_{2n-1}}\, \overline{\s_{2n-2}}$$\cdots$$\overline{\s_1})^{-m}$$T_i$. Thus, we have $\mathcal{L}\cong\mathcal{L}'$.
\end{proof}

\begin{figure}[h]
\begin{center}
\includegraphics[scale=0.7]{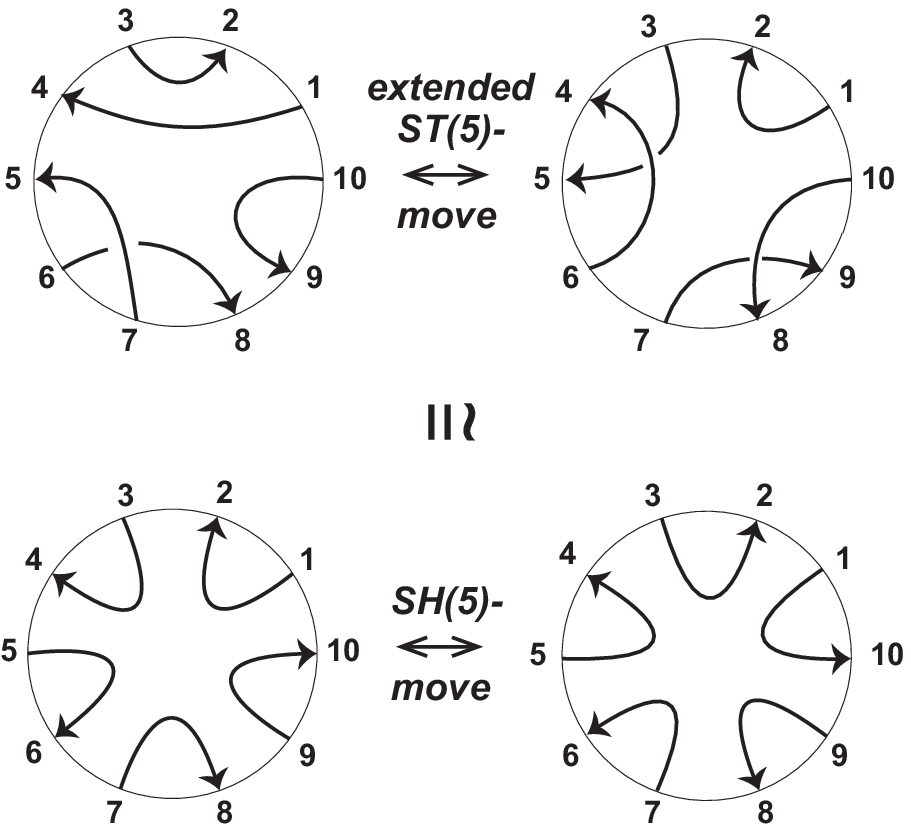}
\end{center}
\caption{}
\label{FIGP5}
\end{figure}

We extend the set of $ST(n)$-moves to the set of local moves that are equivalent to them. A local move that is equivalent to an $ST(n)$-move is called an extended $ST(n)$-move as follows: 
  
\begin{defi} \label{def2-8}
A local move $\mathcal{L}:(D_1, T_1)\leftrightarrow (D_2, T_2)$ is called an extended $ST(n)$-move if there is a local move $\mathcal{T}: (D_3, T_3)\leftrightarrow (D_4, T_4)$ so that $(D_3, T_3)$ and $(D_4, T_4)$ are both trivial, $\mathcal{L}\cong\mathcal{T}$ and $(D_3, T_3)\neq (D_4, T_4)$ (see FIG. \ref{FIGP5}). When we take no notice of the number of arcs in each tangle diagram, an extended $ST(n)$-move is simply called an extended $ST$-move. The sets of extended $ST(n)$ and extended $ST$-moves are denoted by $\mathbb{T}_n$ and $\mathbb{T}$, respectively.
\end{defi}

Clearly an $ST(n)$-move is an extended $ST(n)$-move. 

\begin{defi} \label{def2-9}
Let $\mathcal{L}: (D_1, T_1)\leftrightarrow (D_2, T_2)$ be a local move. Let $p_S$ be a stereographic projection from $D_1$ onto $S^2_{+}\subset S^2$, e.g., $p_S(x,y,0)=\bigl(2x/(1+x^2+y^2),\, 2y/(1+x^2+y^2),\, (1-x^2-y^2)/(1+x^2+y^2)\bigr)$ and $p_N$ is a stereographic projection from $D_2$ onto $S^2_{-}\subset S^2$, e.g., $p_N(x,y,0)=\bigl(2x/(1+x^2+y^2),\, 2y/(1+x^2+y^2),\, (-1+x^2+y^2)/(1+x^2+y^2)\bigr)$. Then $(S^2,\, T_1\cup T_2):=\Bigl( p_S (D_1)\cup p_N (D_2),\, p_S(T_1)\cup p_N(T_2)\Bigr)$ is called the union of the local move $\mathcal{L}$. Here $p_S(T_1)$ and $p_N(T_2)$ are regular diagrams of $T_1$ and $T_2$ into $S^2_{+}$ and $S^2_{-}$, respectively. Each e-point marked $j$ of $(D_1, T_1)$ is identified with an e-point marked $j$ of $(D_2, T_2)$ in the union of $\mathcal{L}$. Each point $\bigl(\cos{ j\over n}\pi , \, \sin { j\over n}\pi ,\, 0\bigr)$ in (the equator of) the union of $\mathcal{L}$ is also marked $j$ and called an e-point in the union of $\mathcal{L}$. The union $(S^2, T_1\cup T_2)$ of $\mathcal{L}$ can be regarded as a link diagram in $S^2$. We call each component of the link diagram ``a component with e-points (in the union of $\mathcal{L}$)." If an e-point in the union of $\mathcal{L}$ is on a component $\mathcal{C}$ with e-points in the union of $\mathcal{L}$, then we call it an e-point of $\mathcal{C}$. Here all e-points in the union of $\mathcal{L}$ are fixed.

\end{defi}

%\begin{defi} \label{def2-10}
%If the union of a local move consists of $m$ components with e-points, we say t%hat the local move has $m$ components' union. In particular, when $m=1$, we say% that the local move has one component's union. 
%\end{defi}

\section{A classification of $ST$-moves}

The $SH(n)$-move, as illustrated in FIG. \ref{shn}, is denoted by $\mathcal{H}_n$. It is easy to see that $\mathcal{H}_n\in\mathbb{T} _n\subset \mathbb{T}\subset \mathbb{L}$. In this section, we shall state some lemmas and Theorem 1. Before describing them, we define some notations.

Let $\mathcal{L}\in \mathbb{L}$ and $\mathcal{C}$ be a component with e-points in the union of $\mathcal{L}$. Let $E(\mathcal{C})$ be the set of integers marked to the e-points of $\mathcal{C}$. We assign the positive integers marked to the e-points of $\mathcal{C}$ to the following notations: Let $A(1)=\min E(\mathcal{C})$. Suppose that the component with an e-point marked $A(1)$ contains $2k$ e-points in total. Now, starting with this e-point, move along the arc in $D_1$. When we first arrive at an e-point, we denote the number marked to the e-point by $A(2)$. Moving on from the e-point marked $A(2)$ along the arc in $T_2$, when we arrive at the next e-point, we denote the positive integer of this e-point by $A(3)$. In this way, we assign the positive integers from $1$ to $2k$ to the notations $A(1), \, A(2), \ldots$ and $A(2k)$. In addition, we define the following notations. Let $B(1)$, $B(2)$, $\cdots$, $B(2k)$ be the integers marked to the e-points of $\mathcal{C}$ satisfying $B(1)<B(2)<\cdots < B(2k)$.

Let $S(\mathcal{C})=\{ \, i \, | \, A(i)-A(i-1)\neq 1, \, i=2,3,\ldots , 2k \}$ and $L(\mathcal{C})=\{\, i\,\, |\,\, B(i)-B(i-1)\neq 1 ,\, i=2,3,\ldots , 2k\}$. If $S(\mathcal{C})\neq\emptyset$ (or $L(\mathcal{C})\neq\emptyset$, resp.), then we let $s(\mathcal{C})=\min S(\mathcal{C})$ (or $l(\mathcal{C})=\min S(\mathcal{C})$, resp.). And then, we apply the following sequence of braiding operations to both tangle diagrams of $\mathcal{L}$:
$$\s_{A(s-1)+1}\s_{A(s-1)+2}\cdots \s_{A(s)-2}\s_{A(s)-1} $$  
$$(\textrm{or}\,\,\, \s_{B(l-1)+1}\s_{B(l-1)+2}\cdots \s_{B(l)-2}\s_{B(l)-1},\,\,\,\textrm{resp.})$$ \\
where $s=s(\mathcal{C})$ and $l=l(\mathcal{C})$. Note that $2\leq s, \,\,  l\leq 2k$, $2\leq A(s)-A(s-1)$ and $2\leq B(l)-B(l-1)$. The local move and the component with e-points obtained from $\mathcal{L}$ and $\mathcal{C}$ using the above sequences of braiding operations, are denoted by $\mathcal{L}'$ and $\mathcal{C}'$, respectively. Let $A'(1), \, A'(2), \ldots ,\, A'(2k)$, $B'(1), \, B'(2), \ldots ,\, B'(2k)$ be the positive integers marked to the e-points of $\mathcal{C}'$ in the manner described above. Let $S(\mathcal{C}')=\{ \, i \, | \, A'(i)-A'(i-1)\neq 1, \, i=2,3,\ldots , 2k \}$ and $L'(\mathcal{C})\,=\left\{\, i\,\, |\,\, B'(i)-B'(i-1)\neq 1 ,\, i=2,3,\ldots , 2k\right\}$. If $S(\mathcal{C}')\neq\emptyset$ (or $L(\mathcal{C}')\neq\emptyset$, resp.), then we let $s(\mathcal{C}')=\min S(\mathcal{C}')$ (or $l(\mathcal{C}')=\min L(\mathcal{C}')$, resp.). Then we see that $s(\mathcal{C})<s(\mathcal{C}')$ (or $l(\mathcal{C})<l(\mathcal{C}')$, resp.).

\begin{prop}\label{prop3-1}
Let $\mathcal{C}$ be a component with $2k$ e-points in the union of a local move $\mathcal{L}$ such that $S(\mathcal{C})\neq\emptyset$. Let $\mathcal{L}'$ and $A'(*)$ be the notations described above. If the following $(a)$ or $(b)$ holds for $\mathcal{L}$, then the following $(a)'$ or $(b)'$ holds for $\mathcal{L}'$:

$(a)$ The tangle diagrams of $\mathcal{L}$ are both trivial.

$(b)$ The condition $(a)$ does not hold. All overcrossings of the tangle diagrams of $\mathcal{L}$ are on the arc, say $\alpha _1$, whose e-points are marked $A(s(\mathcal{C})-1)$ and $A(s(\mathcal{C}))$, and all undercrossings of the diagrams are on arcs except $\alpha _1$. 

$(a)'$ The tangle diagrams of $\mathcal{L}'$ are both trivial.

$(b)'$ The condition $(a)'$ does not hold. Each tangle diagram of $\mathcal{L}'$ is equal to a tangle diagram satisfying the following: All overcrossings are on the arc, say $\beta _2$, whose e-points are marked $A'(s(\mathcal{C}))$ and $A'(s(\mathcal{C})+1)$, and all undercrossing are on arcs except $\beta _2$. 

Here, if $s(\mathcal{C})=2k$, then we let $A'(s(\mathcal{C})+1)):=A'(1)$.

\end{prop}

\begin{figure}[h]
\begin{center}
\includegraphics[scale=0.8]{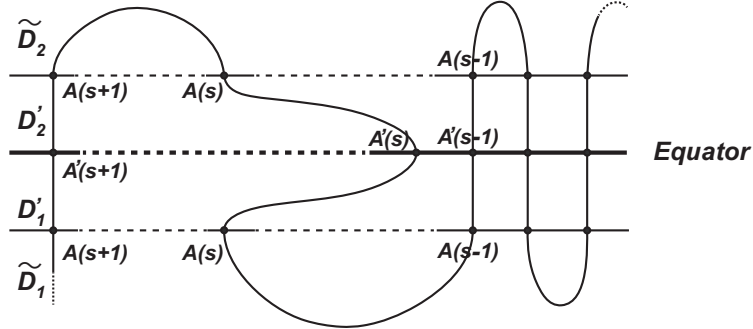}
\end{center}
\caption{Union of $\mathcal{L}'$; Case (I) of $A(s)<A(s+1).$}
\label{prop3-1f}
\end{figure}

\begin{figure}[h]
\begin{center}
\includegraphics[scale=0.8]{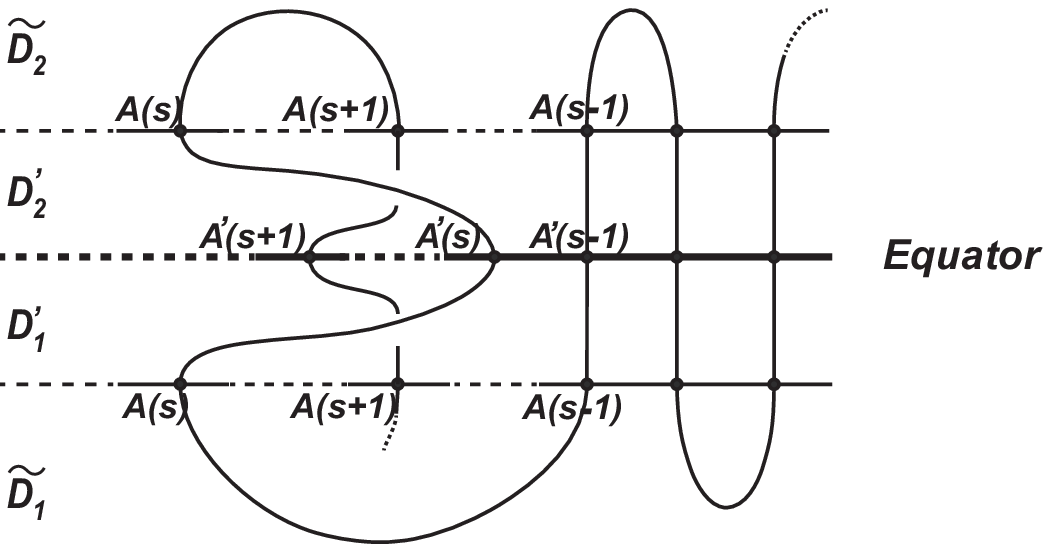}
\end{center}
\caption{Union of $\mathcal{L}'$; Case (I) of $A(s)>A(s+1).$}
\label{prop3-1g}
\end{figure}

\begin{figure}[h]
\begin{center}
\includegraphics[scale=0.7]{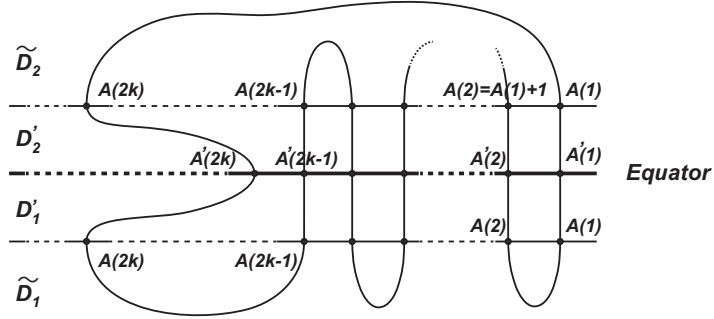}
\end{center}
\caption{Union of $\mathcal{L}'$; Case (II).}
\label{prop3-2}
\end{figure}

\begin{proof}
Let $\mathcal{L}:(D_1, T_1)\leftrightarrow (D_2, T_2)$ and $\mathcal{L}':(D_1', T_1')\leftrightarrow (D_2', T_2')$. Let $(\widetilde{D_i}, \widetilde{T_i})$ be the shrunk diagram of $(D_i, T_i)$ obtained from $(D_i, T_i)$ by the sequence of braiding operations, where $i=1, 2$. We can assume w. l. o. g. that $\alpha_1$ is on $D_1$. Let $\a_1'\in \widetilde{D_1}$ be the arc whose end points are marked $A(s-1)\in\partial\widetilde{D_1}$ and $A(s)\in\partial\widetilde{D_1}$. Let $\a_2\in D_1'\setminus Int(\widetilde{D_1})$ be the arc whose end points are marked $A(s)\in\partial\widetilde{D_1}$ and $A'(s)\in\partial D_1'$, and let $\a_3\in D_2'\setminus Int(\widetilde{D_2})$ be the arc whose end points are marked $A'(s)\in\partial D_2'$ and $A(s)\in\partial\widetilde{D_2}$. Let $\b_1\in D_1'$ be the arc whose end points are marked $A'(s-1)$ and $A'(s)$ (see FIG. \ref{prop3-1f}-\ref{prop3-2}). Here $s=s(\mathcal{C})$.

(I) Case $2\leq s\leq 2k-1$ (see FIG. \ref{prop3-1f} and FIG. \ref{prop3-1g}): Let $\a_4\in \widetilde{D_2}$ be the arc whose end points are marked $A(s)$ and $A(s+1)$. Let $\b_2\in \widetilde{D_2'}$ be the arc whose e-points are marked $A'(s)$ and $A'(s+1)$. 

If the condition $(a)$ holds for $\mathcal{L}$, then $\a_1'$ and $\a_4$ possess no crossings, $\a_2$ and $\a_3$ may do overcrossings and no undercrossings. If the condition $(b)$ holds for $\mathcal{L}$, then $\a_1'$ has overcrossings, $\a_2$ and $\a_3$ may possess overcrossings, $\a_4$ does no crossings. Therefore, even if crossings exist or do not exist on $\a_1'$, all crossings on $\b_1$ are overcrossings and are removed using a finite sequence of Reidemeister moves because the e-points marked $A'(s-1)$ and $A'(s)$ are adjacent. Thus, if a crossing exists in a tangle diagram of $\mathcal{T}'$, then it is on the arc $\a_3$. Therefore, it is on $\b_2$. Remark that if $\b_2$ has a crossing with itself, then only one self-crossing exists and we have $A(s+1)<A(s)$. Since $\a_3$ passes over and $A'(s+1)>A'(s)$ by construction, the only one self-crossing can be removed using Reidemeister moves (see FIG. \ref{prop3-1g}). Hence, the condition $(a)'$ or $(b)'$ holds for $\mathcal{L}'$. 

(II) Case $s=2k$ (see FIG. \ref{prop3-2}): Let $\a_4\in \widetilde{D_2}$ be the arc whose end points are marked $A(s)$ and $A(1)$. Let $\b_2\in \widetilde{D_2'}$ be the arc whose e-points are marked $A'(s)$ and $A(1)$. Even if crossings exist on $\a_1'$, $\a_2$ or $\a_3$, they are removed using a finite sequence of Reidemeister moves because $\b_1$ and $\b_2$ are overpasses, and the e-points marked $A'(2k-1)$ and $A'(2k)$ are adjacent to each other. Thus, the tangle diagrams of $\mathcal{L}'$ are both trivial. Hence, the condition $(a)'$ holds for $\mathcal{L}'$. 
\end{proof}


Note that if $(b)$ or $(b)'$ holds for $\mathcal{T}$ or $\mathcal{T}'$, then $S(\mathcal{C})\neq\emptyset$ or $S(\mathcal{C}')\neq\emptyset$, respectively.

\begin{prop}\label{prop3-2-1}
Let $\mathcal{C}$ be a component with $2k$ e-points in the union of a local move $\mathcal{L}$ such that $L(\mathcal{C})\neq\emptyset$. Let $\mathcal{L}'$ and $B'(*)$ be the notations described in the beginning of this section. If the following $(c)$ or $(d)$ holds for $\mathcal{L}$, then the following $(c)'$ or $(d)'$ holds for $\mathcal{L}'$:

$(c)$ The tangle diagrams of $\mathcal{L}$ are both trivial.

$(d)$ The condition $(c)$ does not hold. All overcrossings of the tangle diagrams of $\mathcal{L}$ are on the arc, say $\alpha _1$, whose e-points are marked $B(l(\mathcal{C})-1)$ and $B(l(\mathcal{C}))$, and all undercrossings of the diagrams are on arcs except $\alpha _1$. 

$(c)'$ The tangle diagrams of $\mathcal{L}'$ are both trivial.

$(d)'$ The condition $(c)'$ does not hold. Each tangle diagram of $\mathcal{L}'$ is equal to a tangle diagram satisfying the following: All overcrossings are on the arc, say $\beta _2$, whose e-points are marked $B'(l(\mathcal{C}))$ and $B'(l(\mathcal{C})+1)$, and all undercrossing are on arcs except $\beta _2$. 

Here, if $l(\mathcal{C})=2k$, then we let $B'(l(\mathcal{C})+1)):=B'(1)$.
\end{prop}

\begin{proof}
This follows from Proposition \ref{prop3-1}.
\end{proof}

Note that if $(d)$ or $(d)'$ holds for $\mathcal{T}$ or $\mathcal{T}'$, then $L(\mathcal{C})\neq\emptyset$ or $L(\mathcal{C}')\neq\emptyset$, respectively.

\begin{lem}\label{lem3-1}
Let $\mathcal{T}\in \mathbb{T}_n$ and $c(\mathcal{T})=1$. Then $\mathcal{T} \cong\mathcal{H}_n$.

\end{lem}

\begin{proof} Let $\mathcal{C}$ be the component with e-points in the union of $\mathcal{T}$. If $\mathcal{T}=\mathcal{H}_n$, then $\mathcal{T} \cong\mathcal{H}_n$. Therefore, suppose that $\mathcal{T}\neq\mathcal{H}_n$, i.e., $\mathcal{T}$ is not an $ST$-move or $S(\mathcal{C})\neq\emptyset$. We may assume that $\mathcal{T}$ is an $ST(n)$-move and $S(\mathcal{C})\neq\emptyset$ because $\mathcal{T}$ is an extended $ST(n)$-move and $\mathcal{T}\neq\mathcal{H}_n$. Since $S(\mathcal{C})\neq 0$ and the condition $(a)$ in Proposition \ref{prop3-1} holds for $\mathcal{T}$, the condition $(a)'$ or $(b)'$ holds for $\mathcal{T}'$ by Proposition \ref{prop3-1}. If $S(\mathcal{C}')=\emptyset$, then the condition $(a)'$ holds for $\mathcal{T}'$. Thus, $\mathcal{T'}=\mathcal{H}_n$. So we have $\mathcal{T} \cong\mathcal{H}_n$. Otherwise, the assumption of Proposition \ref{prop3-1} holds for $\mathcal{T'}$ and that $s(\mathcal{C})<s(\mathcal{C}')$.

Suppose that $S(\mathcal{C}^i)\neq\emptyset$ and the condition $(a)$ or $(b)$ in Proposition \ref{prop3-1} holds for $\mathcal{T}^i$ for $i=1, 2, \ldots , r-1$. Here $\mathcal{T}^0:=\mathcal{T}$, $\mathcal{T}^i:=(\mathcal{T}^{i-1})'$, $\mathcal{C}^0:=\mathcal{C}$ and $\mathcal{C}^i=(\mathcal{C}^{i-1})'$. We let $\mathcal{T}^{r}:=(\mathcal{T}^{r-1})'$ and $\mathcal{C}^r:=(\mathcal{C}^{r-1})'$. Then the condition $(a)'$ or $(b)'$ holds for $\mathcal{T}^r$ by Proposition \ref{prop3-1}. If $S(\mathcal{C}^r)=\emptyset$, then the condition $(a)'$ holds for $\mathcal{T}^r$. Thus, $\mathcal{T}^r=\mathcal{H}_n$. So we have $\mathcal{T}\cong\mathcal{H}_n$. Otherwise, the assumption of Proposition \ref{prop3-1} holds for $\mathcal{T}^r$ and $s(\mathcal{C}^{r-1})<s(\mathcal{C}^r)$. 

Hence, by induction there is a positive integer $m$ such that $S(\mathcal{C}^m)=\emptyset$ since the number of e-points on $\mathcal{C}$ is finite. Then two tangle diagrams of $\mathcal{T}^m$ are trivial. Hence, $\mathcal{T}^m\cong\mathcal{H}_n$ because of $\mathcal{T}\cong\mathcal{T}^1\cong\cdots\cong\mathcal{T}^m\cong\mathcal{H}_n$. This completes the proof.

\end{proof}

\begin{defi} \label{def3-1}
Let $\mathcal{T}\in\mathbb{T}$. Let $\mathcal{C}$ be a component with e-points in the union of $\mathcal{T}$. Let $B(1)$, $B(2)$, $\cdots$, $B(2k)$ be the integers marked to the e-points of $\mathcal{C}$ satisfying $B(1)<B(2)<\cdots < B(2k)$. If $L(\mathcal{C})=\emptyset$, then $\mathcal{C}$ is said to be appropriate. If every component with e-points in the union of $\mathcal{T}$ is appropriate, then $\mathcal{T}$ is said to be appropriate. 
\end{defi}

If the number of components with e-points in the union of $\mathcal{T}$ is one, then $\mathcal{T}$ is appropriate.

\begin{lem} \label{lem3-2} 
Let $\mathcal{T}\in \mathbb{T}$. Then there exists an appropriate $ST$-move that is equivalent to $\mathcal{T}$.
\end{lem}

\begin{proof}
We may suppose that $\mathcal{T}$ is an $ST$-move and is not appropriate. The latter property means that a component $\mathcal{C}$ with e-points in the union of $\mathcal{T}$ exists such that $\mathcal{C}$ is not appropriate, i.e., $L(\mathcal{C})\neq\emptyset$. 

Let $B(1)$, $B(2)$, $\cdots$, $B(2k)$ be the positive integers marked to the e-points of $\mathcal{C}$ satisfying $B(1)<B(2)<\cdots < B(2k)$. Now $L(\mathcal{C})\neq\emptyset$ and the condition $(c)$ in Proposition \ref{prop3-2-1} holds for $\mathcal{T}$. Therefore, $(c)'$ or $(d)'$ in Proposition \ref{prop3-2-1} holds for $\mathcal{T}'$ by Proposition \ref{prop3-2-1}. If $L(\mathcal{C}')=\emptyset$, then the condition $(c)'$ holds for $\mathcal{T}'$. That is, $\mathcal{C}'$ is appropriate and $\mathcal{T}'$ is an $ST$-move.

If $L(\mathcal{C}')\neq \emptyset$, then the condition $(c)'$ or $(d)'$ in Proposition \ref{prop3-2-1} holds for $\mathcal{T}'$. This means that the assumption of Proposition \ref{prop3-2-1} holds for $\mathcal{T}'$. Note that any arc of $\mathcal{C}'$ do not have a crossing with an arc of $\mathcal{C}'$. If a diagram of $\mathcal{T}'$ has a crossing, then its overcrossing is on an arc of $\mathcal{C}'$ and its undercrossing is on an arc of a component except $\mathcal{C}'$. 

Suppose that $L(\mathcal{C}^i)\neq\emptyset$ and the condition $(c)$ or $(d)$ in Proposition \ref{prop3-2-1} holds for $\mathcal{T}^i$ for $i=1, 2, \ldots , r-1$. Here $\mathcal{T}^0:=\mathcal{T}$, $\mathcal{T}^i:=(\mathcal{T}^{i-1})'$, $\mathcal{C}^0:=\mathcal{C}$ and $\mathcal{C}^i:=(\mathcal{C}^{i-1})'$. We let $\mathcal{T}^r:=(\mathcal{T}^{r-1})'$ and $\mathcal{C}^r:=(\mathcal{C}^{r-1})'$. Then the condition $(c)'$ or $(d)'$ holds for $\mathcal{T}^r$ by Proposition \ref{prop3-2-1}. If $L(\mathcal{C}^r)=\emptyset$, then the condition $(c)'$ holds for $\mathcal{T}^r$. Thus, $\mathcal{C}^r$ is appropriate. Otherwise, the assumption of Proposition \ref{prop3-2-1} holds for $\mathcal{T}^r$. Note that any arc of $\mathcal{C}^r$ do not have a crossing with an arc of $\mathcal{C}^r$. If a diagram of $\mathcal{T}^r$ has a crossing, then the overcrossing is on an arc of $\mathcal{C}^r$ and the undercrossing is on an arc of a component except $\mathcal{C}^r$. 

Hence, by induction there is a positive integer $m$ such that $L(\mathcal{C}^m)=\emptyset$ since the number of e-points on $\mathcal{C}$ is finite. Then two tangle diagrams of $\mathcal{T}^m$ are trivial.

If $\mathcal{T}^m$ is not appropriate, then we can continue such an operation until $\mathcal{T}^m$ becomes an appropriate $ST$-move. Additionally, we should remark that every appropriate component with e-points remains appropriate even after the application of these moves.

Hence, there exists an appropriate $ST$-move $\mathcal{T}^m$ that is equivalent to $\mathcal{T}$. This completes the proof.

\end{proof}

Let $\mathcal{L}\in \mathbb{L}$. The number of components with e-points in the union of $\mathcal{L}$ is denoted by $c(\mathcal{L})$. Henceforth, suppose that $j$ is a natural integer such that $1\leq j\leq c(\mathcal{T})$.

\begin{lem} \label{lem3-3}
Let $\mathcal{T}$ be an appropriate $ST(n)$-move. Let $\mathcal{C}_1, \mathcal{C}_2, \ldots , \mathcal{C}_{c(\mathcal{T})}$ be the components with e-points in the union of $\mathcal{T}$. Let $k_j$ be the number of e-points of $\mathcal{C}_j$. Then there exists a finite sequence of braiding operations that transforms $\mathcal{T}$ into an appropriate $ST(n)$-move and $\mathcal{C}_j$ into a component with e-points marked $B_j'(1)$, $B_j'(2)$, $\ldots$, $B_j'(k_j)$ satisfying the following: $1=B_1'(1)<\cdots < B_1'(k_1)< B_2'(1)< \cdots < B_2'(k_2)<\cdots <B_{c(\mathcal{T})}'(1)< \cdots < B_{c(\mathcal{T})}'(k_{c(\mathcal{T})})=2n$. 
\end{lem}

\begin{proof} 
Proposition \ref{prop2-1} ensures the existence of a finite sequence of braiding operations, i.e., a rotation, that transforms the marks of the e-points of $\mathcal{C}_1$ into $1, 2, \ldots$ and $k_1$. The local move and the component with e-points obtained from $\mathcal{T}$ and $\mathcal{C}_j$ using the finite sequence of braiding operations, are once again denoted by $\mathcal{T}$ and $\mathcal{C}_j$, respectively. 

Let $B_j(1)$, $B_j(2)$, $\cdots$, $B_j(k_j)$ be the integers marked to the e-points of $\mathcal{C}_j$ satisfying $B_j(1)<B_j(2)<\cdots < B_j(k_j)$. Then $\mathcal{T}$ is an appropriate $ST(n)$-move, and the following property holds: $1=B_1(1)<2=B_1(2)<\cdots k_1=B_1(k_1)$. Let $B(1)=B_1(1)$, $B(2)=B_1(2)$, $\ldots$, $B(k_1)=B_1(k_1)$, $B(k_1+1)=B_2(1)$, $\ldots$, $B(k_1+k_2)=B_2(k_2)$ and
$$L(\mathcal{C}_1\cup \mathcal{C}_2)\,=\left\{\, i\,\, |\,\, B(i)-B(i-1)\neq 1 , \, i=2,3,\ldots , k_1+k_2\right\}.$$ 
Now we see that $L(\mathcal{C}_j)=\emptyset$ for each $j$.

By the proof of Lemma \ref{lem3-2}, there exists a finite sequence of braiding operations that transforms $\mathcal{C}_j$ into $\mathcal{C}_j'$ such that $L(\mathcal{C}_1'\cup \mathcal{C}_2')=\emptyset$. Then the integers $B_j'(1)<B_j'(2)<\cdots <B_j'(k_i)$ marked to the e-points of $\mathcal{C}_j'$ satisfy the following property: $1=B_1'(1)<\cdots < B_1'(k_1)< B_2'(1)< \cdots < B_2'(k_2)=k_1+k_2$. Then we see that the local move $\mathcal{T}'$ obtained from $\mathcal{T}$ is an $ST$-move. 

Suppose that $\mathcal{T}$ is an appropriate $ST(n)$-move and $1=B_1(1)< B_2(1)< \cdots <B_{M-1}(1)<2n$, where $2\leq M\leq c(\mathcal{T})$. Let $B(1)=B_1(1)$, $B(2)=B_1(2)$, $\ldots$, $B(k_1)=B_1(k_1)$, $B(k_1+1)=B_2(1)$, $\ldots$, $B(k_1+k_2)=B_2(k_2)$, \ldots , $B(k_1+k_2+\cdots +k_{M-1}+1)=B_M(1)$, $\ldots$ , $B(k_1+k_2+\cdots +k_M)=B_M(k_M)$ and
$$L(\mathcal{C}_1\cup \mathcal{C}_2\cup \cdots \cup \mathcal{C}_M)\,=\left\{\, i\,\, |\,\, B(i)-B(i-1)\neq 1 , \, i=2,3,\ldots , k_1+k_2+\cdots k_M\right\}.$$

By the proof of Lemma \ref{lem3-2}, there exists a finite sequence of braiding operations that transforms $\mathcal{C}_j$ into $\mathcal{C}_j'$ such that $L(\mathcal{C}_1'\cup \mathcal{C}_2'\cup \cdots \cup \mathcal{C}_M')=\emptyset$. Then $\mathcal{T}'$ is an appropriate $ST$-move, and the integers $B_j'(1)<B_j'(2)<\cdots <B_j'(k_i)$ marked to the e-points of $\mathcal{C}_j'$ satisfy the following property: $1=B_1'(1)<\cdots < B_1'(k_1)< B_2'(1)< \cdots < B_2'(k_2)<\cdots <B_{M}'(1)< \cdots < B_M'(k_M)=k_1 +k_2+\cdots k_M$. Hence, by induction the proof is completed.

\end{proof}

\begin{defi} \label{def3-2}
Let $\mathcal{T}$ be an $ST$-move. Let $\mathcal{C}$ be a component with e-points in the union of $\mathcal{T}$. Let $A(1)$, $A(2), \ldots ,\, A(2k)$ be the integers marked to the e-points of $\mathcal{C}$ as described in the beginning of this section. If $A(1)$ is an initial e-point of $\mathcal{C}$ and $S(\mathcal{C})=\emptyset$, then we say that this component $\mathcal{C}$ is an $SH(+)$-type. If the components with e-points in the union of $\mathcal{T}$ are all $SH(+)$-types, then $\mathcal{T}$ is said to be an $SH(+)$-type.
\end{defi} 

In Definition \ref{def3-2}, if $\mathcal{C}$ is an $SH(+)$-type, then its form is uniquely decided and we have $A(i)=B(i)$ for $i=1,2,\ldots , 2k$. Here the notation $B(i)$ is the integers marked to the e-points of $\mathcal{C}$ as in Definition \ref{def3-1}. 

\begin{lem}\label{lem3-4}
Let $\mathcal{T}$ be an appropriate $ST$-move. Then there is an $SH(+)$-type that is equivalent to $\mathcal{T}$. 
\end{lem}

\begin{proof}
Suppose that $\mathcal{T}$ is not an $SH(+)$-type. That is, there is a component, say $\mathcal{C}$, with e-points in the union of $\mathcal{T}$ such that $\mathcal{C}$ is not an $SH(+)$-type. Choose an initial e-point of $C$. By Proposition \ref{prop2-1}, there exists a finite sequence of braiding operations that transforms the e-point into the e-point marked $1$. By the proof of Lemma \ref{lem3-1}, the component $C$ with e-points can be transformed into an $SH(+)$-type using a finite sequence of braiding operations without affecting the other components with e-points in the union of $\mathcal{T}$. 

We can continue such an operation until every component with e-points in the union becomes an $SH(+)$-type. We complete the proof. 
\end{proof}

\begin{figure}[h]
\begin{center}
\includegraphics[scale=0.7]{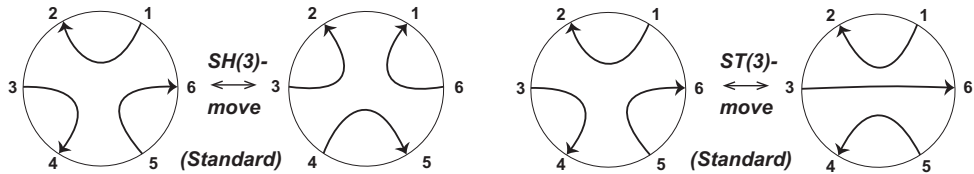}
\end{center}
\caption{Standard $ST(3)$-moves.}
\label{st3m}
\end{figure}

\begin{defi} \label{def3-3}
Let $\mathcal{T}$ be an $SH(+)$-type. Let $\mathcal{C}_1$, $\mathcal{C}_2$, $\cdots $, $\mathcal{C}_{c(\mathcal{T})}$ be the components with e-points in the union of $\mathcal{T}$.  Let $A_j(1)$, $A_j(2)$, $\cdots$, $A_j(k_j)$ be the integers marked to the e-points of $\mathcal{C}_j$ satisfying $A_1(1)< A_2(1)< \cdots < A_{c(\mathcal{T})}(1)$. If $k_1\leq k_2\leq\cdots\leq k_{c(\mathcal{T})}$, then $\mathcal{T}$ is said to be standard (see FIG. \ref{st3m}).
\end{defi}

\begin{defi}\label{def3-4}
Let $\mathcal{T}\in \mathbb{T}_n.$ Let $\mathcal{C}_1$, $\mathcal{C}_2$, $\cdots $, $\mathcal{C}_{c(\mathcal{T})}$ be the components with e-points in the union of $\mathcal{T}$ and $n_j={1\over 2}|E(\mathcal{C}_j)|$. If $n_1\leq n_2\leq \cdots\leq n_{c(\mathcal{T})}$, then the partition $n=n_1+n_2+\cdots +n_{c(\mathcal{T})}$ of $n$ is called the arc-decomposition of $\mathcal{T}$, where $|\ast |$ is the number of the set $\ast$.
\end{defi}

By the definition, we see that a standard $ST(n)$-move $\mathcal{T}$, whose arc-decomposition is $n=n_1+n_2+\cdots +n_{c(\mathcal{T})}$, is uniquely decided. 

The next theorem states that there is a one-to-one correspondence between the set of equivalence classes of extended $ST(n)$-moves and the set of standard $ST(n)$-moves for each $n$. 

\begin{thm} \label{thm1}
Given $\mathcal{T}\in\mathbb{T}_n$, there is a uniquely standard $ST(n)$-move that is equivalent to $\mathcal{T}$. Conversely, given a standard $ST(n)$-move, there is an extended $ST(n)$-move that is equivalent to the standard one.
\end{thm}

\begin{proof}
Let $\mathcal{T}\in\mathbb{T}_n$. Let $\mathcal{C}_1, \mathcal{C}_2, \ldots , \mathcal{C}_{c(\mathcal{T})}$ be the components with e-points in the union of $\mathcal{T}$, in which the number of e-points of $\mathcal{C}_{i+1}$ is greater than or equal to one of $\mathcal{C}_i$, where $i=1, 2, \ldots , c(\mathcal{T})-1$. Then there is a standard $ST(n)$-move $\mathcal{T}'$ such that $\mathcal{T}'\cong\mathcal{T}$, by Lemmas \ref{lem3-2}, \ref{lem3-3} and \ref{lem3-4}. Suppose that there is a standard $ST(n)$-move $\mathcal{T}''$ such that $\mathcal{T}''\neq\mathcal{T}'$ and $\mathcal{T}''\cong\mathcal{T}$. Then the arc-decomposition of $\mathcal{T}''$ is not equal to the arc-decomposition of $\mathcal{T}'$, which is a contradiction. Because a braiding operation on a local move $\mathcal{L}$ does not change the number $c(\mathcal{L})$ and the number of arcs in each tangle diagram of any component with e-points of $\mathcal{L}$. Thus, $\mathcal{T}'$ is uniquely decided. 

Conversely, given a standard $ST(n)$-move, it is also an extended $ST(n)$-move. Thus, Theorem 1 holds. 
\end{proof}

\begin{figure}[h]
\begin{center}
\includegraphics[scale=0.7]{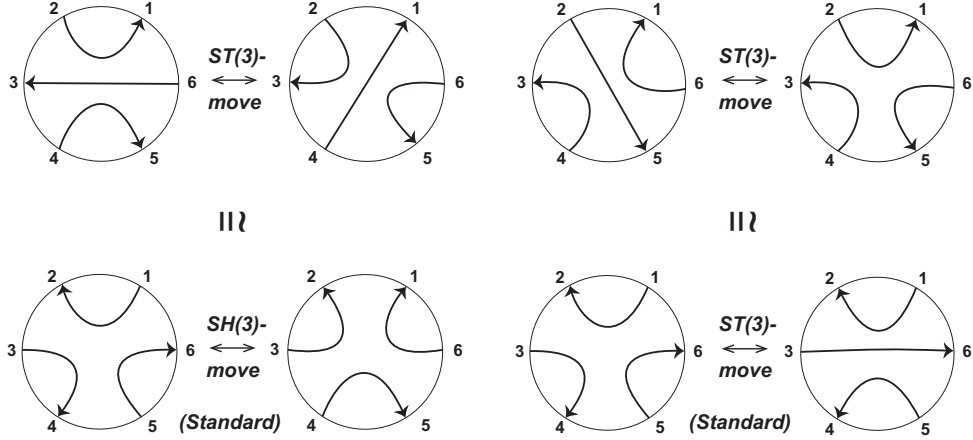}
\end{center}
\caption{Equivalent moves}
\label{st3n}
\end{figure}

For example, for $n=3$, there exist two standard $ST(3)$-moves in all, i.e., the two local moves in the second row of FIG. \ref{st3n}, which are not equivalent. Although there are many extended $ST(3)$-moves, each of them is equivalent to the diagram on the lower left or the diagram on the lower right in FIG. \ref{st3n}. A left-side local move and a right-side local move of FIG. \ref{st3n} are not equivalent. On the other hand, the two $ST(3)$-moves, on the left or on the right of FIG. \ref{st3n} are equivalent. Thus, the four $ST(3)$-moves as in FIG. \ref{st3n} are classified into two equivalence classes whose representatives are the lower-side local moves of FIG. \ref{st3n}. 

By Theorem \ref{thm1}, we can choose a standard $ST$-move as a representative of an equivalence class of extended $ST$-moves. The standard $ST(n)$-move whose arec-decomposition is $n=n_1+n_2+\cdots + n_r$ is denoted by $\langle n_1, n_2, \cdots , n_r\rangle$, where $n_1\leq n_2\leq\cdots\leq n_r$. 

Next, we have a necessary and sufficient condition for two extended $ST$-moves to be equivalent.

\begin{cor} \label{cor1} Let $\mathcal{T}_n\in \mathbb{T}_n$ and $\mathcal{T}_m\in \mathbb{T}_m$. Let $n=n_1+n_2+\cdots + n_r$ and $m=m_1+m_2+\cdots + m_s$ be the arc-decompositions of $\mathcal{T}_n$ and $\mathcal{T}_m$, respectively. Then $\mathcal{T}_n\cong\mathcal{T}_m$ if and only if $r=s$ and $n_i=m_i$, where $i=1,2,\ldots , r$.
\end{cor}

\begin{proof} Suppose that $\mathcal{T}_n\cong\mathcal{T}_m$. Since a braiding operation on a local move $\mathcal{L}$ does not change the number $c(\mathcal{L})$ and the number of arcs in each tangle diagram of any component with e-points of $\mathcal{L}$, the equalities $r=s$ and $n_i=m_i$ hold, where $i=1,2,\ldots , r$. 

Suppose that $r=s$ and $m_i=n_i$ holds, where $i=1,2,\ldots , r$. By Theorem 1, there is a uniquely standard $ST(n)$-move $\mathcal{T}$, whose arc-decomposition is $n=n_1+n_2+\cdots n_r$, such that $\mathcal{T}\cong\mathcal{T}_n$ and $\mathcal{T}\cong\mathcal{T}_m$. Therefore, we have $\mathcal{T}_n\cong\mathcal{T}_m$. Hence, we complete the proof. 
\end{proof}

The next corollary provides the number of equivalence classes of extended $ST(n)$-moves. 

\begin{cor} \label{cor2} Let $(2\leq)$ $n$ be a positive integer and $p(n)$ be the partition number of $n$. Then we have $|[\mathbb{T}_n]|=p(n)-1$, where $[*]$ is the set of equivalence classes of the set $*$.
\end{cor}

\begin{proof} This is an immediate consequence of Theorem 1. 
\end{proof}

For example, the number of equivalence classes of extended $ST(3)$-moves is $2$ because the partition number of $3$ is $3$.

\section{A partial order on equivalence classes}

In this section, we shall define a specific operation, called a connecting operation, that transforms an $n$-tangle diagram into an $(n-1)$-tangle diagram in order to describe a binary relation on the set of local moves.

\begin{figure}[h]
\begin{center}
\includegraphics[scale=0.6]{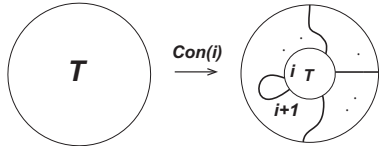}
\end{center}
\caption{Connecting operation $Con(i)$.}
\label{connect-oper}
\end{figure}

\begin{defi} \label{def4-1}
Let $i\, (\leqq 2n)$ be a positive integer and $A:=\{(x,y,0)\in\mathbb{R}^3 \,|\, 1/2\leqq x^2+y^2\leqq 1\}$. Let $\Hat{U}_i$ be a union of $2n-1$ pairwise disjoint arcs $u_k$ embedded properly in $B\setminus Int(\widetilde{B})$, satisfying the following:

\noindent $(I)$ Case of $i=1,2,\ldots , 2n-1.$ 
\begin{itemize}
\item For $j\in\{1, 2, \ldots , i-1\}$, $u_j$ is an arc connecting the points $1/2(\cos{ j\over n}\pi , \, \sin {j\over n}\pi,$ $0)$ and $(\cos{ j\over n-1}\pi , \, \sin { j\over n-1}\pi, \, 0)$ whose image $p(u_j)$ is on the annulus $A$. 

\item For $j\in\{i+2, i+3, \ldots , 2n\}$, $u_j$ is an arc connecting the points $1/2(\cos{ j\over n}\pi,$ $\sin {j\over n}\pi , \, 0)$ and $(\cos{ j-2\over n-1}\pi , \, \sin { j-2\over n-1}\pi , \, 0)$ whose image $p(u_j)$ is on $A$.

\item $u_i$ is an arc connecting the points $1/2(\cos{ i\over n}\pi , \, \sin { i\over n}\pi , \, 0)$ and $1/2(\cos{i+1\over n}\pi ,$ $\sin {i+1\over n}\pi$, $\, 0)$ whose image $p(u_i)$ is on $A$.

\item The diagram $(A, U_i:=p(\Hat{U}_i))$ has no crossings. 
\end{itemize}

\noindent $(II)$ Case of $i=2n.$ 
\begin{itemize}
\item For $j\in\{2,3,\ldots , 2n-1\}$, $u_j$ is an arc connecting the points $1/2(\cos{ j\over n}\pi,$ $\sin {j\over n}\pi , \, 0)$ and $(\cos{ j-1\over n-1}\pi , \, \sin { j-1\over n-1}\pi , \, 0)$ whose image $p(u_j)$ is on $A$.

\item $u_{2n}$ is an arc connecting the points $1/2(\cos{1\over n}\pi , \, \sin {1\over n}\pi , \, 0)$ and $1/2(1, \, 0 , \, 0)$ whose image $p(u_{2n})$ is on $A$.

\item The diagram $(A, U_{2n}:=p(\Hat{U}_{2n}))$ has no crossings. 
\end{itemize}
If $(D, T)$ is an $n$-tangle diagram, the e-point marked $i$ is an initial e-point $($or a terminal e-point$)$, the e-point marked $i+1$ $($$1$ if $i=2n$$)$ is a terminal e-point $($or an initial e-point$)$, and $(\widetilde{D}, \widetilde{T})$ is the shrunk diagram of $(D, T)$, then we can regard $(D', T'):=(A\cup\widetilde{D},\, U_i\cup\widetilde{T})$ as an $(n-1)$-tangle diagram. The operation that transforms $(D, T)$ into $(D', T')$ is called the connecting operation $Con(i)$ on $(D,T)$ and we write $T'=Con(i)\, T$. $($see FIG. \ref{connect-oper}$)$. Here the orientation of $(D', T')$ is induced from the orientation of $(D, T)$.

\end{defi}

Next, we define a binary relation $\preceq$ on the set of local moves.

\begin{defi} \label{def4-2}
Let $\mathcal{L}: T_1\leftrightarrow T_2$ and $\mathcal{L}': T_1'\leftrightarrow T_2'$ be local moves. If there is a finite sequence of connecting and/or braiding operations $o_1$, $o_2$, $\cdots$, $o_r$ such that $T_1'=o_r\cdots o_2o_1T_1$ and $T_2'=o_r\cdots o_2o_1T_2$, then the operation that transforms $\mathcal{L}$ into $\mathcal{L}'$ is called a sequence of connecting and/or braiding operations $o_1$, $o_2$, $\cdots$, $o_r$ on $\mathcal{L}$ and we say that $\mathcal{L}'$ can be obtained (or realized) from $\mathcal{L}$ by a finite sequence of connecting and/or braiding operations $o_1$, $o_2$, $\cdots$, $o_r$, and we write $\mathcal{L}\preceq\mathcal{L}'$. In particular, if $o_1$, $o_2$, $\cdots$, $o_r$ are all braiding operations, then we say that $\mathcal{L}'$ can be obtained (or realized) from $\mathcal{L}$ by a finite sequence of braiding operations $o_1$, $o_2$, $\cdots$, $o_r$, and we write $\mathcal{L}'\cong\mathcal{L}$ or $\mathcal{L}'=o_r\cdots o_2o_1\mathcal{L}$. 

Let $\mathcal{L}\in\mathbb{L}_n$. If there is a positive integer $k_1$ such that $o_i=\overline{\sigma_{k_i}}$, $o_r=Con(k_r)$, $i=1,2,\ldots , r-1$, 
where
$$ k_{i+1}\begin{cases}%
=k_i +1 \qquad\qquad\qquad \textrm{if}\,\,\,\,  k_i+1\leqq 2n,  \\
\equiv k_i +1 \pmod{2n} \quad \textrm{if}\,\,\,\,  k_i+1 >2n,  
\end{cases} $$
then we call the sequence of operations $o_1$, $o_2$, $\cdots$, $o_r$, a sequence of connecting and braiding operations on $\mathcal{L}$ that connects the e-points marked $k_1$ and $k_r+1$.
\end{defi}

\begin{figure}[h]
\begin{center}
\includegraphics[scale=0.5]{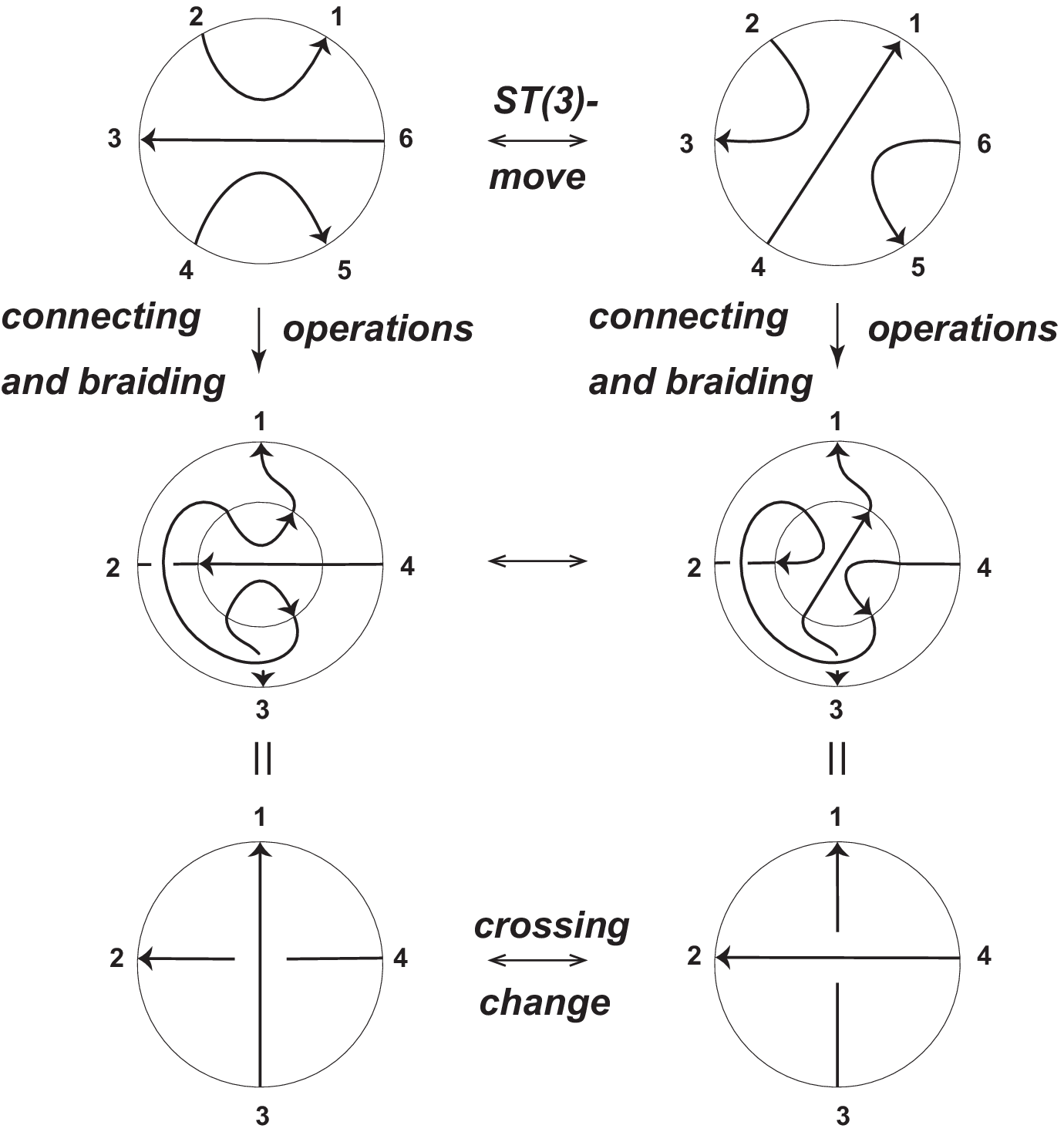}
\end{center}
\caption{\,}
\label{real}
\end{figure}

\begin{figure}[h]
\begin{center}
\includegraphics[scale=0.5]{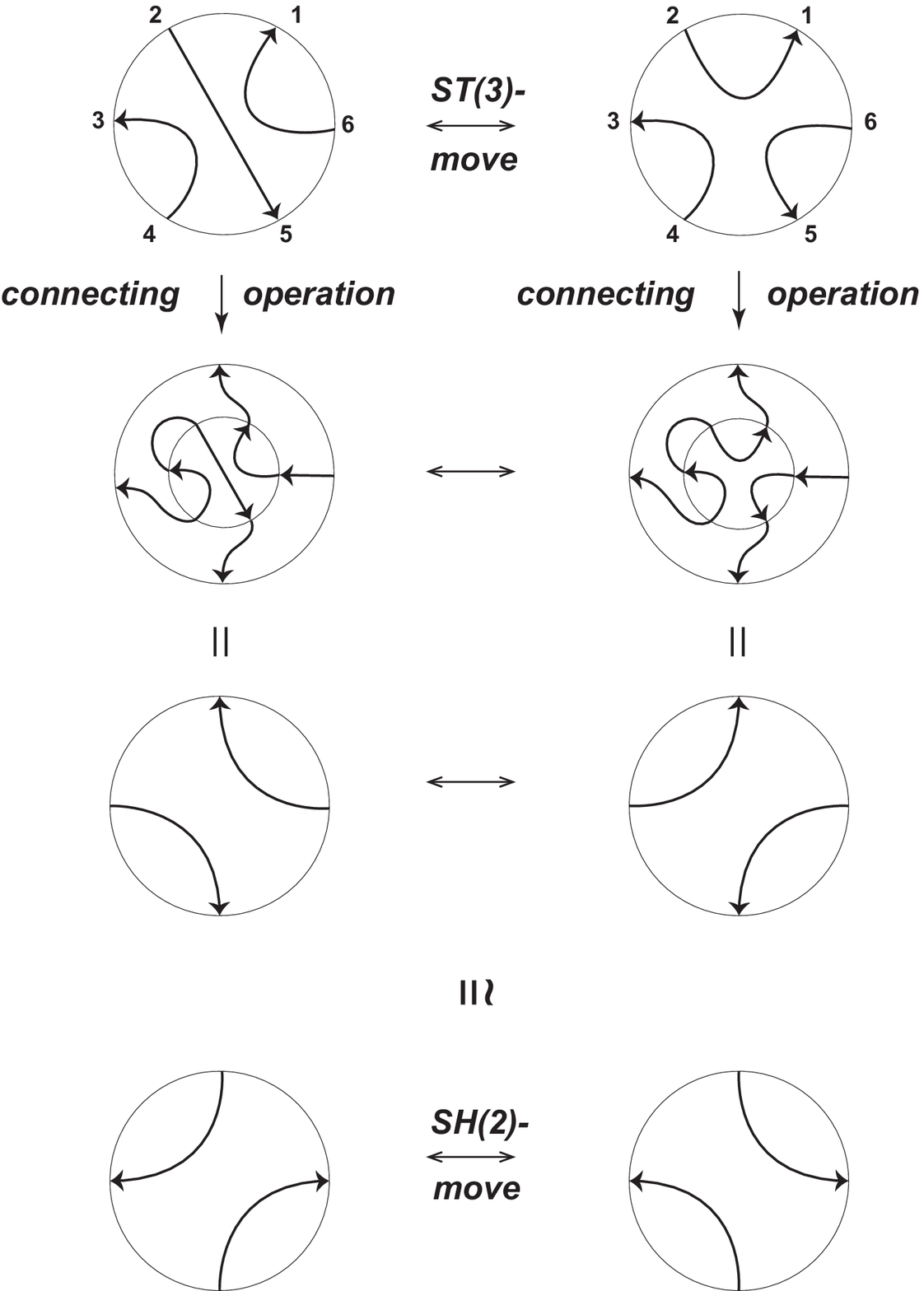}
\end{center}
\caption{\,}
\label{real2}
\end{figure}

Let $\mathcal{X}$ be the $X$-move (crossing change), $\mathcal{H}_2$ be the $SH(2)$-move, $\mathcal{T}$ be the $ST(3)$-move as shown in the upper-side diagram of FIG. \ref{real} and $\mathcal{T}'$ be the $ST(3)$-move as in the upper-side diagram of FIG. \ref{real2}. Then we see that $\mathcal{T}\preceq\mathcal{X}$ and $\mathcal{T}'\preceq\mathcal{H}_2$, which can be seen in FIG. \ref{real} and FIG. \ref{real2}, respectively. Because $\mathcal{X}=Con(4)\, \overline{\sigma_3}\, \overline{\sigma_2}\, \mathcal{T}$ and $\mathcal{H}_2=Con(2)\mathcal{T}'$.


The $SH(n)$-move, as illustrated in the diagram in FIG. \ref{shn} is denoted by $\mathcal{H}_n$. The sets of extended $ST(n)$-moves, extended $ST$-moves and local moves are denoted by $\mathbb{T}_n$, $\mathbb{T}$ and $\mathbb{L}$, respectively. The subset of $\mathbb{L}$ whose element is a pair of $n$-tangle diagrams is denoted by $\mathbb{L}_n$. 

Let $[\mathcal{L}]:=\{\mathcal{L}'\in\mathbb{L}\,|\, \mathcal{L}\cong\mathcal{L}'\,\}$ denote the equivalence class to which $\mathcal{L}$ belongs. Then the relation $\preceq$ is well-defined under the equivalence relation $\cong$ as follows: Suppose that $\mathcal{L}_1\cong\mathcal{L}'_1$ and $\mathcal{L}_2\cong\mathcal{L}'_2$. If $\mathcal{L}_1\preceq\mathcal{L}_2$, then we only need to prove that $\mathcal{L}'_1\preceq\mathcal{L}'_2$. Since $\mathcal{L}_1\cong\mathcal{L}'_1$ and $\mathcal{L}_2\cong\mathcal{L}'_2$, there are finite sequences of braiding operations $\s_{i_1}, \s_{i_2}, \cdots , \s_{i_r}$ and $\s_{j_1}, \s_{j_2}, \cdots , \s_{j_n}$ such that $\mathcal{L}_1=\s_{i_r}\cdots \s_{i_2}\s_{i_1}\mathcal{L}'_1$ and $\mathcal{L}_2'=\s_{i_r}\cdots \s_{i_2}\s_{i_1}\mathcal{L}_2$. Therefore we see that $\mathcal{L}_1'\preceq \s_{i_r}\cdots \s_{i_2}\s_{i_1}\mathcal{L}'_1=\mathcal{L}_1\preceq\mathcal{L}_2\preceq\s_{i_r}\cdots \s_{i_2}\s_{i_1}\mathcal{L}_2=\mathcal{L}_2'$. Hence, $\mathcal{L}'_1\preceq\mathcal{L}'_2$.

Next, we show that the binary relation $\preceq$ is a partial order on the set $[\mathbb{L}]:=\{\, [\mathcal{L}]\,|\, \mathcal{L}\in\mathbb{L}\,\}$ of equivalence classes of local moves. Namely, the relation $\preceq$ on the set $[\mathbb{L}]$ is reflexive, antisymmetric and transitive. This means it satisfies the following $(i)-(iii)$ for any $\mathcal{F}$, $\mathcal{G}$ and $\mathcal{L}$ in $\mathbb{L}$: $(i)$ $[\mathcal{F}]\preceq[\mathcal{F}]$, $(ii)$ if $[\mathcal{F}]\preceq[\mathcal{G}]$ and $[\mathcal{G}]\preceq[\mathcal{L}]$, then $[\mathcal{F}]\preceq[\mathcal{L}]$, and $(iii)$ if $[\mathcal{F}]\preceq[\mathcal{G}]$ and $[\mathcal{G}]\preceq[\mathcal{F}]$, then $[\mathcal{F}]=[\mathcal{G}]$. According to Definition \ref{def4-1}, we see that $(i)$ and $(ii)$ hold. Let $\mathcal{F}\in\mathbb{L}_n$ and $\mathcal{G}\in\mathbb{L}_m$. Suppose that $[\mathcal{F}]\preceq[\mathcal{G}]$ and $[\mathcal{G}]\preceq[\mathcal{F}]$. It holds that $m\leqq n$ by the relation $\mathcal{F}\preceq\mathcal{G}$ and it holds that $n\leqq m$ by the relation $\mathcal{G}\preceq\mathcal{F}$. Therefore, we have $n=m$. Namely, $\mathcal{G}$ (or $\mathcal{F}$, resp.) can be obtained from $\mathcal{F}$ (or $\mathcal{G}$, resp.) by a finite sequence of braiding operations. Therefore, we have $\mathcal{F}\cong\mathcal{G}$. Thus, we have $[\mathcal{F}]=[\mathcal{G}]$ and so $(iii)$ holds. Hence, $(\mathbb{L}, \preceq)$ is a partially ordered set. We will regard each standard $ST(n)$-move as a representative of an equivalence class of an extended $ST(n)$-move. 

\section{some results}

In this section, we discuss partially orderings between equivalence classes of extended $ST$-moves. We have already the partially orderings between equivalence classes of $SH$-moves. The following result can be found in \cite{HNT}.

\begin{lem}\label{lem5-1}$($\cite[Lemma 16]{HNT}$)$
For any $n$, we have $\mathcal{H}_{2n+1}\preceq \mathcal{H}_{2n-1}$ and $\mathcal{H}_{2n}\preceq \mathcal{H}_{2n-1}$.
\end{lem}

\begin{figure}[h]
\begin{center}
\includegraphics[scale=0.7]{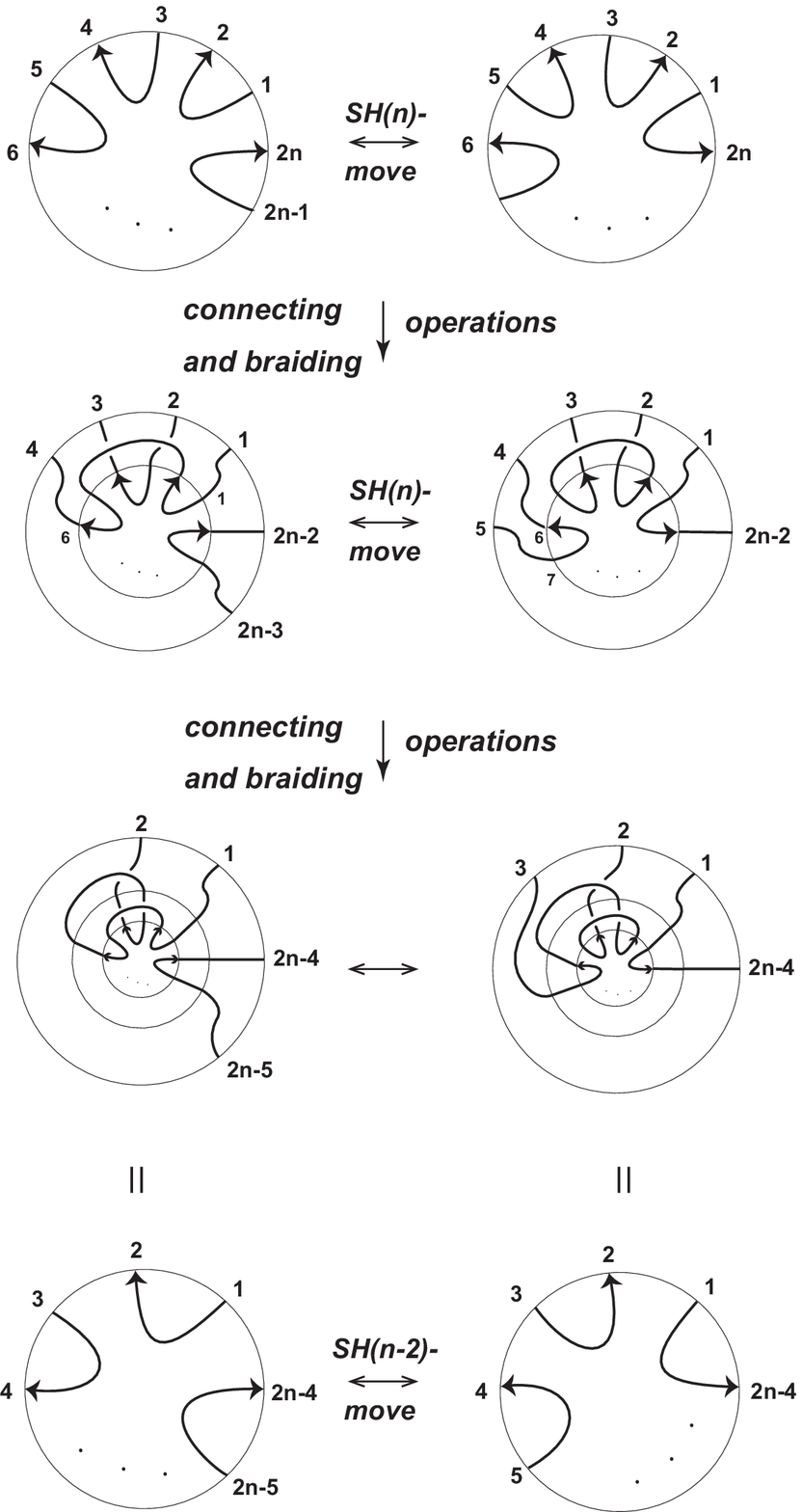}
\end{center}
\caption{\,}
\label{SH1}
\end{figure}

\begin{proof}
Using the sequence of two connecting operations shown in FIG. \ref{SH1}, the $SH(n-2)$-move can be realized from the $SH(n)$-move. 
\end{proof}

Therefore, we have the following corollary.

\begin{cor}\label{cor3}
For any $n$, we have $[\mathcal{H}_{2n+1}]\preceq [\mathcal{H}_{2n-1}]$ and $[\mathcal{H}_{2n}]\preceq [\mathcal{H}_{2n-2}]$.
\end{cor}

Note that if $\mathcal{L}\in\mathbb{L}_n$, then for any $\mathcal{L}'\in [\mathcal{L}]$, $\mathcal{L}'\in\mathbb{L}_n$ and $c(\mathcal{L}')=c(\mathcal{L})$.

\begin{prop}\label{prop5-1}
Let $\mathcal{L}_n\in\mathbb{L}_n$ and $\mathcal{L}_m\in\mathbb{L}_m.$ If $[\mathcal{L} _n]\preceq[\mathcal{L}_m]$, then $m-c(\mathcal{L}_m)\equiv n-c(\mathcal{L}_n)\pmod{2}$ and $m\leqq n$.
\end{prop}

\begin{proof} We assume that $[\mathcal{L} _n]\preceq [\mathcal{L}_m]$ i.e. $\mathcal{L}_m$ can be realized from $\mathcal{L}_n$ by a finite sequence of connecting and/or braiding operations. A braiding operation on $\mathcal{L}_n$ does not change the number $n-c(\mathcal{L}_n)$.  A connecting operation changes $n$ and $c(\mathcal{L}_n)$ into $n-1$ and $c(\mathcal{L}_n)\pm 1$, respectively. Therefore, $m\leqq n$ and it changes $(n-c(\mathcal{L}_n))$ into $(n-c(\mathcal{L}_n))$ or $(n-c(\mathcal{L}_n)-2)$. Thus, $m-c(\mathcal{L}_m)=n-c(\mathcal{L}_n)-2l$ for some $l\in\{0, 1, 2, \cdots\}$, and hence we have $m-c(\mathcal{L}_m)\equiv n-c(\mathcal{L}_n)\pmod{2}.$ We complete the proof.
\end{proof}

The next corollary gives a necessary and sufficient conditions for the relation between equivalence classes of $SH$-move to exist. 

\begin{cor}\label{cor4}
$[\mathcal{H} _m]\preceq[\mathcal{H}_n]$ if and only if $n\equiv m \pmod{2}$ and $n\leqq m$.
\end{cor}

\begin{proof}
This follows from Proposition \ref{prop5-1} and Lemma \ref{lem5-1}.

\end{proof}

By Corollary \ref{cor4}, $\Bigl(\cup _{n=1}^\infty [\mathcal{H}_{2n+1}], \,\preceq\Bigr)$ and $\Bigl(\cup_{n=1}^\infty [\mathcal{H}_{2n}], \preceq\Bigr)$ are both totally ordered sets. Also, we can conclude from Proposition \ref{prop5-1} that there are no binary relations $\preceq$ between an element of $\cup _{n=1}^\infty\mathcal{H}_{2n}$ and an element of $\cup _{n=1}^\infty\mathcal{H}_{2n+1}$. 

\begin{cor}\label{cor5}
Let $\mathcal{T} _m\in\mathbb{T} _m$, $\mathcal{T}_n\in\mathbb{T}_n$ and $c(\mathcal{T}_m)=c(\mathcal{T}_n)=1$. Then $[\mathcal{T} _m]\preceq [\mathcal{T}_n]$ if and only if $n\equiv m \pmod{2}$ and $n\leqq m$.
\end{cor}

\begin{proof}
This follows from Corollaries \ref{cor1} and \ref{cor4}.
\end{proof}

The following two lemmas show necessary and sufficient conditions for the relation between the equivalence class of an $ST$-move and one of an $SH$-move to exist.

\begin{lem}\label{lem5-2}
Let $\mathcal{T} _m\in\mathbb{T} _m$, $\mathcal{T} _n\in\mathbb{T} _n$ and $c(\mathcal{T}_n)=1$. Then $[\mathcal{T} _m]\preceq \, [\mathcal{T} _n]$ if and only if $n\equiv m-c(\mathcal{T} _m)+1\pmod{2}$ and $n\leqq m-c(\mathcal{T} _m)+1$.
\end{lem}

\begin{proof} When $c(\mathcal{T}_m)=1$, the proposition holds by Corollary \ref{cor5}. Therefore, we only need to prove the case in which $c(\mathcal{T}_m)\neq 1$.

Suppose that $[\mathcal{T} _m]\preceq \, [\mathcal{T}_n]$, i.e. $\mathcal{T} _n$ can be realized from $\mathcal{T}_m$ by using connecting and/or braiding operations. Since the number of components in the union of $\mathcal{T}_m$ must be changed into one to construct $\mathcal{T}_n$ from $\mathcal{T}_m$ by connecting and/or braiding operations, at least $c(\mathcal{T}_m)-1$ connecting operations must be applied on $\mathcal{T}_m$. If $c(\mathcal{T}_m)-1$ connecting operations are applied on $\mathcal{T}_m$, then the number of arcs of two tangle diagrams of new $\mathcal{T}_m$ must be $m-\{c(\mathcal{T}_m)-1\}$. Therefore, by Corollary \ref{cor5}, it is necessary for $\mathcal{T}_n$ to be realized from $\mathcal{T}_m$ that $n\equiv m-c(\mathcal{T} _m)+1\pmod{2}$ and $n\leqq m-c(\mathcal{T}_m)+1.$

Suppose that $n\equiv m-c(\mathcal{T}_m)+1\pmod{2}$ and $n\leqq m-c(\mathcal{T}_m)+1$. Let $\mathcal{H}_n=\langle n\rangle\in [\mathcal{T}_n]$ and $\mathcal{T}_m'=\langle k_1, k_2, \ldots , k_{c(\mathcal{T}_m)}\rangle\in [\mathcal{T}_m]$ be the standard $ST$-moves. Then we show that  $\mathcal{H}_n \preceq \mathcal{T}_m'$. Since $n-1\leqq m-c(\mathcal{T}_m)$, the inequality $n<m$ holds. Thus, $\mathcal{H}_n$ may be realized from $\mathcal{T}_m'$ by a finite sequence of connecting and/or braiding operations.

A sequence of connecting operations must transform $c(\mathcal{T}_m)$ components with e-points of $\mathcal{T}_m'$ into a single component with e-points. In order to connect them, we need at least $(c(\mathcal{T}_m)-1)$ times of connecting operations on $\mathcal{T}_m'$. Let $N=m-\{c(\mathcal{T}_m)-1\}$.

(i) Case $n=N$: Apply the following finite sequence of connecting operations $Con(k_1)$, $Con(k_1+k_2)$, $\cdots$, $Con(k_1+k_2+\cdots k_{c(\mathcal{T}_m)-2})$, $Con(k_1+k_2+\cdots k_{c(\mathcal{T}_m)-1})$ on $\mathcal{T}'_m$. Then $\mathcal{T}_m'$ can be transformed into $\mathcal{H}_n$. Thus, we have $\mathcal{T} _m'\preceq \, \mathcal{H}_n$. 

(ii) Case $n< N$ and $n\equiv N\pmod{2}$: Then there is a positive integer $l$ such that $n=N-2l$. Therefore, we have $\mathcal{H}_n=\mathcal{H}_{N-2l}$. Further, the following relation can be obtained from Case (i) and Lemma \ref{lem5-1}: $\mathcal{T}' _m\preceq \, \mathcal{H} _N\preceq \, \mathcal{H} _{N-2}\preceq \cdots \preceq\, \mathcal{H} _{N-2l}=\mathcal{H}_n$. Thus, we have $\mathcal{T}' _m\preceq \, \mathcal{H}_n$. 

Hence, a necessary and sufficient condition is for $\mathcal{T}_n$ to be realized from $\mathcal{T}_m$ that $n\equiv m-c(\mathcal{T} _m)+1\pmod{2}$ and $n\leqq m-c(\mathcal{T} _m)+1.$ This completes the proof. 
\end{proof}

\begin{lem}\label{lem5-3}
Let $\mathcal{T} _m\in\mathbb{T} _m$, $\mathcal{T} _n\in\mathbb{T} _n$ and $c(\mathcal{T}_n)=1$. Then $[\mathcal{T} _n]\preceq \, [\mathcal{T} _m]$ if and only if $m\equiv n-c(\mathcal{T} _m)+1\pmod{2}$ and $m\leqq n-c(\mathcal{T} _m)+1$. 
\end{lem}

\begin{proof} When $c(\mathcal{T}_m)=1$, the proposition holds by Corollary \ref{cor5}. Therefore, we assume that $c(\mathcal{T}_m)\neq 1$. 

Suppose that $[\mathcal{T} _n]\preceq \, [\mathcal{T} _m]$, i.e. $\mathcal{T} _m$ can be realized from $\mathcal{T}_n$. Then we see that $m-c(\mathcal{T} _m)\equiv n-1\pmod{2}$ by Proposition \ref{prop5-1}. Since the number of components in the union of $\mathcal{T}_n$ must be changed into $c(\mathcal{T}_m)$ by the sequence of connecting operations, at least $c(\mathcal{T}_m)-1$ connecting operations must be applied on $\mathcal{T}_n$. If $c(\mathcal{T}_m)-1$ connecting operations are applied on $\mathcal{T}_n$, then the number of arcs of two tangle diagrams of new $\mathcal{T}_n$ must be $n-\{c(\mathcal{T}_m)-1\}$. Therefore, we have $m\leqq n-\{c(\mathcal{T}_m)-1\}$.

Suppose that $m\equiv n-c(\mathcal{T}_m)+1\pmod{2}$ and $m\leqq n-c(\mathcal{T} _m)+1$. Let $\mathcal{H}_n=\langle n\rangle\in [\mathcal{T}_n]$ and $\mathcal{T}_m'=\langle k_1, k_2, \ldots , k_{c(\mathcal{T}_m)}\rangle\in [\mathcal{T}_m]$ be the standard $ST$-moves. Then we show that $\mathcal{H}_n \preceq \mathcal{T}_m'$. Since $m-1\leqq n-c(\mathcal{T}_m)$, the inequality $m<n$ holds. Thus, $\mathcal{T}'_m$ may be realized from $\mathcal{H}_n$ by using connecting operations.

A sequence of connecting operations can transform a single component with e-points $\mathcal{H}_n$ into $c(\mathcal{T}'_m)=c(\mathcal{T}_m)$ components with e-points. In order to divide the single component of $\mathcal{H}_n$ into $c(\mathcal{T}'_m)$ components, at least $(c(\mathcal{T}'_m)-1)$ times of connecting operations must be applied on $\mathcal{H}_n$. Let $N=m+c(\mathcal{T}'_m)-1$.

(i) Case $n=N$: Apply sequences of connecting and braiding operations on $\mathcal{H}_n$ that connect the e-points marked $k_1+1$ and $2n$, $k_1+1$ and $k_1+k_2+2$, $k_1+k_2+1$ and $k_1+k_2+k_3+2$, $k_1+k_2+k_3+1$ and $k_1+k_2+k_3+k_4+2$, $\cdots$ and $k_1+k_2+\cdots +k_{c(\mathcal{T}_m)-2}+1$ and $k_1+k_2+\cdots +k_{c(\mathcal{T}_m)-1}+2$. Thus, $\mathcal{H}_n$ can be transformed into $\mathcal{T}'_m$ i.e. $\mathcal{H}_n\preceq \, \mathcal{T}' _m$. 

(ii) Case $N<n$ and $N\equiv n\pmod{2}$: Then there is a positive integer $l$ such that $N=n-2l$. Therefore, we have $\mathcal{H}_n=\mathcal{H}_{N+2l}$. Further, the following relation can be obtained from Lemma \ref{lem5-1}: $\mathcal{H}_n=\mathcal{H}_{N+2l}\preceq \cdots \preceq\, \mathcal{H} _{N}$. Using Case $(i)$, we have $\mathcal{H}_N\preceq \, \mathcal{T}' _m$. Thus, $\mathcal{H}_n\preceq \, \mathcal{T}' _m$.

Hence, a necessary and sufficient condition is for $\mathcal{T}_m$ to be realized from $\mathcal{T}_n$ that $n-\{c(\mathcal{T}_m)-1\}\equiv m\pmod{2}$ and $m\leqq n-\{c(\mathcal{T}_m)-1\}$ i.e. $m+c(\mathcal{T} _m)-1\leqq n.$ This completes the proof. 
\end{proof}

From Proposition \ref{prop5-1}, we see that any finite sequence of connecting or braiding operations on a local move $\mathcal{L}$ does not change the number $n-c(\mathcal{L})$ modulo $Z_2$. Therefore, we can define as follows.

\begin{defi}\label{X-type}
Let $\mathcal{T}\in\mathbb{T}_n$. If the integer $n-c(\mathcal{T})$ is even, then we say that $\mathcal{T}$ is an $X$-type. If $n-c(\mathcal{T})$ is odd, then we say that $\mathcal{T}$ is an $O$-type.
\end{defi}

\begin{thm}\label{thm2}
Let $\mathcal{T}\in\mathbb{T}_n$. If $\mathcal{T}$ is an $X$-type, then $\mathcal{T}$ can realize the ordinary unknotting operation. Otherwise, $\mathcal{T}$ can realize the $SH(2)$-move. 
\end{thm}

\begin{proof}

We note that $1\leqq n-c(\mathcal{T})\leqq n-1$ because of $1\leqq c(\mathcal{T})\leqq n-1$. Therefore, if $n-c(\mathcal{T})$ is an even integer i.e. $\mathcal{T}$ is an $X$-type, then $3\leqq n-c(\mathcal{T})+1\equiv 1\pmod{2}$. Lemma \ref{lem5-2} gives us $\mathcal{T}\preceq\mathcal{H}_{n-c(\mathcal{T})+1}.$ Therefore, from \cite{HNT}, we see that $\mathcal{H}_{n-c(\mathcal{T})+1}\preceq\mathcal{H}_3\preceq\mathcal{X}$ and so we have $\mathcal{T}\preceq\mathcal{X}$. 

If $n-c(\mathcal{T})$ is an odd integer i.e. $\mathcal{T}$ is an $O$-type, then $2\leqq n-c(\mathcal{T})+1\equiv 0\pmod{2}$. Lemma \ref{lem5-2} tells us that $\mathcal{T}\preceq\mathcal{H}_{n-c(\mathcal{T})+1}.$ Therefore, from \cite{HNT}, we see that $\mathcal{H}_{n-c(\mathcal{T})+1}\preceq\mathcal{H}_2$. Thus, we have $\mathcal{T}\preceq\mathcal{H}_2$. This completes the proof.

\end{proof}

%\begin{cor}\label{cor5}
%Any extended $ST$-move in $\mathbb{T}_X$ is an unknotting operation.
%\end{cor}

%\begin{proof}
%This corollary follows from Theorem \ref{thm2}.

%\end{proof}

We define $\mathbb{T}_X:=\{\mathcal{T}\in\mathbb{T}\, |\, \mathcal{T}\,\, \textrm{is an $X$-type }\}$ and $\mathbb{T}_O:=\{\mathcal{T}\in\mathbb{T}\, |\, \mathcal{T} \,\, \textrm{is an $O$-type }\}$.

\begin{thm}\label{thm3}
Any local move that realizes an extended $ST$-move in $\mathbb{T}_X$ is an unknotting operation.
\end{thm}

\begin{proof}
This theorem follows from Theorem \ref{thm2}.

\end{proof}

\section{unknotting numbers of $ST$-moves}

Let $K$ be an oriented knot in the 3-sphere $S^3$. If an $ST$-move is an $O$-type, then the number of components of $K$ must change when the $ST$-move is applied to a diagram of $K$. So we treat only $X$-type $ST$-moves.

Let $\mathcal{T}\in\mathbb{T}_X$. We denote the minimum number of $\mathcal{T}$ that can transform a diagram of $K$ into a trivial knot diagram by $u_{\mathcal{T}}(K)$, where the minimum is taken over all diagrams of $K$. The proceeding properties follow from section 3. 

\begin{remark} Let $\mathcal{T}, \mathcal{T'}\in\mathbb{T}_X$. If $\mathcal{T}\preceq\mathcal{T'}$, then we have $u_{\mathcal{T}}(K)\leqq u_{\mathcal{T'}}(K)$ for any oriented knot $K$. Hence, if $\mathcal{T}\cong\mathcal{T'}$, then we have $u_{\mathcal{T}}(K)=u_{\mathcal{T'}}(K)$ for any oriented knot $K$. Because if $\mathcal{T}\cong\mathcal{T'}$, then we have $\mathcal{T}\preceq\mathcal{T'}$ and $\mathcal{T'}\preceq\mathcal{T}$.
\end{remark}

\begin{thm} \label{thm4} Let $\mathcal{T}\in [\langle s_1, s_2, \ldots , s_n\rangle]\subset \mathbb{T}_X, \mathcal{T}'\in [\langle s'_1, s'_2, \ldots , s'_m\rangle]\subset\mathbb{T}_X$. If $\sum_{i=1}^n(s_i-1)\leqq\sum_{i=1}^m(s'_i-1)$, then we have $u_{\mathcal{T}}(K)\leqq u_{\mathcal{T'}}(K)$ for any oriented knot $K$. In particular, if $\sum_{i=1}^n(s_i-1)=\sum_{i=1}^m(s'_i-1)$, then we have $u_{\mathcal{T}}(K)=u_{\mathcal{T'}}(K)$ for any oriented knot $K$.
\end{thm}

\begin{proof}
Let $K$ be an oriented knot, $N=\sum_{i=1}^n(s_i-1)$ and $M=\sum_{i=1}^n(s'_i-1)$. From Lemma \ref{lem5-2}, we see that $\mathcal{T}\preceq \mathcal{H}_{N+1}$. Also Remark tells us that $u_{\mathcal{T}}(K)\leqq u_{\mathcal{H}_{N+1}}(K)$ and $u_{\mathcal{T}'}(K)\leqq u_{\mathcal{H}_{M+1}}(K)$. 

On the other hands, if $u_{\mathcal{T}}(K)=l$, then we have $u_{\mathcal{H}_{N+1}}(K)\leqq l$. Because we can gather one root of each band near one point of the trivial knot and $l$ times of $SH(N+1)$-moves can produce the trivial knot from $K$, we see that $u_{\mathcal{H}_{N+1}}(K)\leqq u_{\mathcal{T}}(K)$. Hence, $u_{\mathcal{T}}(K)= u_{\mathcal{H}_{N+1}}(K)$. Similarly, we have $u_{\mathcal{T}'}(K)= u_{\mathcal{H}_{M+1}}(K)$. 

Therefore, if $N\leqq M$, then we have $u_{\mathcal{T}'}(K)=u_{\mathcal{H}_{M+1}}(K)\leqq u_{\mathcal{H}_{N+1}}(K)=u_{\mathcal{T}'}(K)$ by Lemma \ref{lem5-1} and Remark. Thus, we have $u_{\mathcal{T}}(K)\leqq u_{\mathcal{T}'}(K)$. In particular, if $N=M$, then we have $u_{\mathcal{T}'}(K)=u_{\mathcal{H}_{N+1}}(K)=u_{\mathcal{T}'}(K)$. Thus, we have $u_{\mathcal{T}}(K)=u_{\mathcal{T}'}(K)$. We complete the proof.
 
\end{proof}

In the next proposition, we show that for any oriented knot $K$, there is an $ST$-move so that $K$ can be transformed into the trivial knot by the single $ST$-move.

\begin{thm}\label{thm5} Let $\mathcal{T}\in [\langle s_1, s_2, \ldots , s_n\rangle]\subset \mathbb{T}_X, \mathcal{T}'\in [\langle s'_1, s'_2, \ldots , s'_m\rangle]\subset\mathbb{T}_X$. If $\sum_{i=1}^n(s_i-1)=u_{\mathcal{T}'}(K)\cdot\sum_{i=1}^m(s'_i-1)$, then we have $u_{\mathcal{T}}(K)=1$ for any oriented knot $K$.
\end{thm}

\begin{proof}
Let $K$ be an oriented knot and $N=u_{\mathcal{T}'}(K)\cdot\sum_{i=1}^m(s'_i-1)$. Then $u_{\mathcal{H}_{N+1}}(K)=1$. Because we can gather one root of each band near one point of the trivial knot as in the proof of Theorem \ref{thm4}, the single $SH(N+1)$-move can produce the trivial knot from $K$. If $\sum_{i=1}^n(s_i-1)=N$, then Theorem 3 tells us that $u_{\mathcal{T}}(K)=u_{\mathcal{H}_{N+1}}(K)$. Thus, we have $u_{\mathcal{T}}(K)=1$.

\end{proof}

Let $\mathcal{T}'\in [\langle s'_1, s'_2, \ldots , s'_m\rangle ]\subset\mathbb{T}_X$ and $K$ be an oriented knot. From the definition of $ST$-moves we see that the value $\sum_{i=1}^m(s'_i-1)$ is a natural number. Therefore, the value $u_{\mathcal{T}'}(K)\cdot\sum_{i=1}^m(s'_i-1)$ is also a natural number. Let $p$ be any natural number. Then there exist natural numbers $s_1$ and $s_2$ such that $s_1+s_2-2=p$. Thus, there is an $ST$-move (e.g. $(s_1, s_2)$, $s_1+s_2-2=u_{\mathcal{T}'}(K)\cdot\sum_{i=1}^m(s'_i-1)$) except $SH$-moves so that $K$ can be transformed into the trivial knot by the single $ST$-move.

\section{some examples}

\begin{figure}[h]
\begin{center}
\includegraphics[scale=0.7]{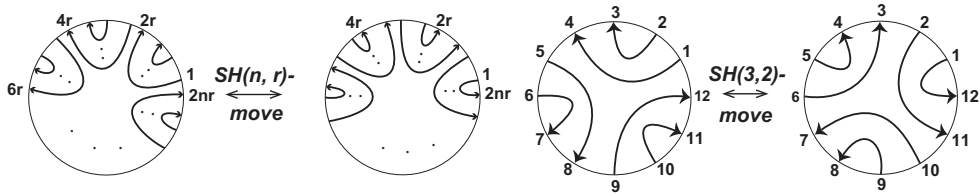}
\end{center}
\caption{$SH(n,r)$-move.}
\label{fig4}
\end{figure}

Let $\mathbb{N}$ be the set of natural numbers and let $r\in \mathbb{N}$. The local move, shown in FIG. \ref{fig4}, is called an $SH(m,n)$-move and denoted by $\mathcal{H}(n,r)$. In particular, $\mathcal{H}(n,1)=\mathcal{H} _n$. 

From Corollary \ref{cor1}, we have the following. 
           
\begin{example} $\mathcal{H}(n,r)\preceq \, \mathcal{H}(n,s)$ if and only if the following conditions $(i)$ or $(ii)$ holds:
$(i)$ $n$ is even, $s\equiv r \pmod 2$ and $s\leqq r,$
$(ii)$ $n$ is odd and $s\leqq r$.
\end{example}

\noindent Additionally, Lemmas \ref{lem5-2} and \ref{lem5-3} give us the following example.
\begin{example}
For any positive integers $a$ and $b$ $(2\leqq a)$, we have $\mathcal{H}_{ab+b-1}\preceq \, \mathcal{H}(a,b)\preceq \, \mathcal{H}_{ab-b+1}$.

\end{example}
\noindent Lastly, the following example results from Theorem 3.
\begin{example}
Let $a, c\in\mathbb{N}$ and $2\leqq b, d\in\mathbb{N}$. If $b(a-1)=d(c-1)$, then we have $u_{\mathcal{H}(a,b)}(K)=u_{\mathcal{H}(c,d)}(K)=u_{\mathcal{H}_{ab-b+1}}(K)$ for any oriented knot $K$.

\end{example}

\section*{Acknowledgements}
The author would like to thank Professor Y. Nakanishi and Professor S. Satoh for their helpful discussions and suggestions and Professor S. Fukuhara and Professor H. A. Miyazawa for reading a draft of the paper and providing useful comments.

\section*{Acknowledgements}
The author would like to thank Professor H. A. Miyazawa for helpful discussions and suggestions and Professor S. Fukuhara for reading a draft of the paper and for providing useful comments.


\begin{thebibliography}{99}

\bibitem{ai} H. Aida, The oriented $\Delta\sb{ij}$-moves on links, Kobe J. Math. {\bf 9}\ (1992), 163--170.

\bibitem{miya} H. Aida, Unknotting operations of polygonal type, Tokyo J. Math. {\bf 15}\ (1992), 111--121.

\bibitem{co} J. Conway, An Enumeration of Knots and Links, and Some of Their Algebraic Properties, Computational Problems in Abstract Algebra. Oxford, England. Pergamon Press\ (1970), 329--358. 

\bibitem{HNT}
J. Hoste, Y. Nakanishi and K. Taniyama, Unknotting Operation Involving Trivial Tangles, Osaka J. Math. {\bf 27}\ (1990), 555--566.

\bibitem{kawauchi} A. Kawauchi, A survey of knot theory, Birkh{\"a}user Verlag, Basel, 1996.


\bibitem{mura} H. Murakami, Some metrics on classical knots, Math. Ann. {\bf 270}\ (1985), 35--45.

\bibitem{MN}
H. Murakami and N. Nakanishi, On a certain move generating link-homology, Math. Ann. {\bf 284}\ (1989), 75--89.

\bibitem{na}
M. Nagura, Unknotting Operations by Using Oriented Trivial Tangle Diagrams, J. Knot Theory and Its Ramifications {\bf 8} (1999), no. 7, 901--929.

\bibitem{N}
Y. Nakanishi, Replacements in the Conway third identity, Tokyo J. Math. {\bf 14}\ (1991), 197--203. 

\bibitem{re} K. Reidemeister, Elementare Begrundung der Knotentheorie, Abh. Math. Sem. Univ. Hamburg {\bf 5} (1926), 24--32 

\end{thebibliography}
\end{document}